\documentclass[11pt,twoside]{amsart}

\usepackage[english]{babel}
\usepackage{csquotes}
\usepackage{amsfonts}
\usepackage{hyperref}
\usepackage{etoolbox}
\usepackage{amsmath}
\usepackage{amssymb}
\usepackage{amsthm}
\usepackage{stmaryrd}
\usepackage{dsfont}
\usepackage{pifont}
\usepackage{caption}
\usepackage{subcaption}
\usepackage{mathrsfs}
\usepackage{graphicx}
\usepackage{enumerate}
\usepackage[all]{xy}
\usepackage{enumitem}
\usepackage{geometry}
\usepackage{amsfonts}
\usepackage{fancyhdr}
\usepackage{calc}
\usepackage{cancel}
\usepackage{accents}
\usepackage{abraces}
\usepackage{amsaddr}
\usepackage{orcidlink}
\usepackage[new]{old-arrows}
\usepackage{tikz,tikz-3dplot}
\usetikzlibrary{decorations.markings}
\usetikzlibrary{shapes,arrows,plotmarks,patterns,calc,decorations.markings,quotes,angles,math}
\tikzset{
    set arrow inside/.code={\pgfqkeys{/tikz/arrow inside}{#1}},
    set arrow inside={end/.initial=>, opt/.initial=},
    /pgf/decoration/Mark/.style={
        mark/.expanded=at position #1 with
        {
            \noexpand\arrow[\pgfkeysvalueof{/tikz/arrow inside/opt}]{\pgfkeysvalueof{/tikz/arrow inside/end}}
        }
    },
    arrow inside/.style 2 args={
        set arrow inside={#1},
        postaction={
            decorate,decoration={
                markings,Mark/.list={#2}
            }
        }
    },
}
\usepackage[backend=biber,style=alphabetic,sorting=nty,maxnames=50]{biblatex}

\usepackage{color}
\definecolor{linkColor}{rgb}{0.0,0.0,0.554}
\definecolor{citeColor}{rgb}{0.0,0.0,0.554}
\definecolor{fileColor}{rgb}{0.0,0.0,0.554}
\definecolor{urlColor}{rgb}{0.0,0.0,0.554}
\definecolor{promptColor}{rgb}{0.0,0.0,0.589}
\definecolor{brkpromptColor}{rgb}{0.589,0.0,0.0}
\definecolor{gapinputColor}{rgb}{0.589,0.0,0.0}
\definecolor{gapoutputColor}{rgb}{0.0,0.0,0.0}
\usepackage{fancyvrb}
\usepackage{adjustbox}
\usepackage{helvet}

\definecolor{cof}{RGB}{219,144,71}
\definecolor{pur}{RGB}{186,146,162}
\definecolor{greeo}{RGB}{91,173,69}
\definecolor{greet}{RGB}{52,111,72}

\tdplotsetmaincoords{70}{165}

\newcommand{\changefont}{%
    \fontsize{8}{8}\selectfont
}
\mathchardef\mhyphen="2D

\title[Dirichlet--Voronoi domain and injectivity radius of flag manifolds]{Dirichlet--Voronoi domain and injectivity radius of flag manifolds - equivariant cell structure on $O(3)/O(1)^3$}

\author{Arthur Garnier \orcidlink{0000-0003-4069-3203}}
\address{LAMFA, Universit\'e de Picardie Jules Verne, CNRS UMR 7352, \\33, rue Saint-Leu, 80000, Amiens, France.}
\email{arthur.garnier@math.cnrs.fr}

\theoremstyle{plain}
\newtheorem{prop}{Proposition}[section]
\newtheorem{prop-def}[prop]{Proposition-Definition}

\newtheorem{lem}[prop]{Lemma}
\newtheorem{theo}[prop]{Theorem}
\newtheorem{cor}[prop]{Corollary}
\newtheorem{rem}[prop]{Remark}

\newtheorem{definition}[prop]{Definition}
\newtheorem{exemple}[prop]{Example}

\newtheorem*{prop*}{Proposition}
\newtheorem*{prop-def*}{Proposition-Definition}
\newtheorem*{propri*}{Property}
\newtheorem*{lem*}{Lemma}
\newtheorem*{theo*}{Theorem}
\newtheorem*{cor*}{Corollary}
\newtheorem*{rem*}{Remark}
\newtheorem*{rems*}{Remarks}
\newtheorem*{definition*}{Definition}
\newtheorem*{exemple*}{Example}
\newtheorem*{notation*}{Notation}
\newtheorem*{conjecture*}{Conjecture}

\newcommand{\lra}{\longrightarrow}

\newcommand{\ra}{\rightarrow}
\newcommand{\sdp}{\times\kern-.2em\vrule height1.1ex depth-.05ex}
\newcommand{\epi}{\lra \kern-.8em\ra}
\newcommand{\C}{{\mathbb C}}

\newcommand{\N}{{\mathbb N}}
\newcommand{\R}{{\mathbb R}}

\newcommand{\Z}{{\mathbb Z}}

\newcommand{\ho}{\mathrm{Hom}\,}

\newcommand{\longto}{\longrightarrow}

\DeclareMathOperator{\Exp}{\mathrm{Exp}}

\newcommand{\Sph}{\mathbb{S}}

\newcommand{\Pro}{\mathbb{P}}
\newcommand{\rk}{\mathrm{rk}\,}

\newcommand{\tr}{\mathrm{tr}\,}

\newcommand{\Sym}{\mathfrak{S}}

\newcommand{\interior}[1]{\accentset{\circ}{#1}}

\DeclareMathOperator\leng{L}
\DeclareMathOperator\inj{inj}
\DeclareMathOperator\diam{diam}
\DeclareMathOperator\pr{pr}

\DeclareRobustCommand\longtwoheadrightarrow
     {\relbar\joinrel\twoheadrightarrow}

\newcommand{\eq}[1][r]
{\ar@<-3pt>@{-}[#1]
\ar@<-1pt>@{}[#1]|<{}="gauche"
\ar@<+0pt>@{}[#1]|-{}="milieu"
\ar@<+1pt>@{}[#1]|>{}="droite"
\ar@/^2pt/@{-}"gauche";"milieu"
\ar@/_2pt/@{-}"milieu";"droite"}

 \normalbaselines
\makeatletter
\newlength\@SizeOfCirc%
\newcommand{\CircleArrowRight}[1]{%
    \setlength{\@SizeOfCirc}{\maxof{\widthof{#1}}{\heightof{#1}}}%
    \tikz [x=1.0ex,y=1.0ex,line width=.12ex]%
        \draw [->,anchor=center]%
            node (0,0) {#1}%
            (0,0.8\@SizeOfCirc) arc (85:-240:0.8\@SizeOfCirc);%
}%
\newcommand{\CircleArrowLeft}[1]{%
    \setlength{\@SizeOfCirc}{\maxof{\widthof{#1}}{\heightof{#1}}}%
    \tikz [x=1.0ex,y=1.0ex,line width=.12ex]%
        \draw [<-,anchor=center]%
            node (0,0) {#1}%
            (0,0.8\@SizeOfCirc) arc (85:-240:0.8\@SizeOfCirc);%
}%
\makeatother
\makeatletter
\newcommand{\opnorm}{\@ifstar\@opnorms\@opnorm}
\newcommand{\@opnorms}[1]{%
  \left|\mkern-1.5mu\left|\mkern-1.5mu\left|
   #1
  \right|\mkern-1.5mu\right|\mkern-1.5mu\right|
}
\newcommand{\@opnorm}[2][]{%
  \mathopen{#1|\mkern-1.5mu#1|\mkern-1.5mu#1|}
  #2
  \mathclose{#1|\mkern-1.5mu#1|\mkern-1.5mu#1|}
}
\makeatother


\date{\today}

\begin{document}

\begin{abstract}
In the first part of this work, we study Dirichlet--Voronoi domains for discrete isometry groups of Riemannian manifolds, in view of constructing cell structures on homogeneous (complete, real) flag manifolds, equivariant with respect to the action of the Weyl group. We give general results, allowing us to build such a structure from an admissible one on the domain. In particular, the injectivity radius plays a key role in the method.

The second part starts with the computation of the injectivity radius of (real and complex) flag manifolds; a first step towards the application of the method developed in the first part. Then, with the help of the quaternion algebra, we investigate the particular case of the flag manifold $O(3)/O(1)^3$ of $SL_3(\R)$: we prove that the results of the first part apply and derive a new $\Sym_3$-equivariant cell structure on it, whose cellular complex of $\Z[\Sym_3]$-modules is determined.
\end{abstract}

\maketitle


\section*{Introduction}


Our general aim is to construct cell structures for flag manifolds which are equivariant with respect to the Weyl group. In the present paper, we propose an approach using tools from Riemannian geometry. It has the potential to ultimately give a general proof, but several obstacles still need to be overcome. Here we tackle one of them: computing the \textit{injectivity radius} of the flag manifold. The answer turns out to be nice and closely related to the combinatorics of the corresponding root system; see the discussion below and Theorem \ref{injrad} for a precise statement.

Specifically, let $K$ be a simple compact Lie group, and $T<K$ be a maximal torus. Since $T$ is abelian, the \textit{Weyl group} $W=N_K(T)/T$ acts freely on the \textit{flag manifold} $\mathcal{F}:=K/T$, a smooth manifold which carries many structures. Namely, recall that the \textit{complexification} $G=K^\C$ of $K$ is a complex algebraic group that contains $K$ as a (maximal) compact subgroup. We also say that $K$ is the \textit{compact real form} of $G$. Choosing a Borel subgroup $B<G$ containing $T$, the \textit{Iwasawa decomposition} yields a diffeomorphism $\mathcal{F}\simeq G/B$, so that $\mathcal{F}$ also has the structure of a smooth complex projective variety. Moreover, as a complex variety, $\mathcal{F}$ has a preferred real structure. Indeed, the complex group $G$ also possesses a \textit{split real form} $G_\R$, thus giving a real structure on $G$ and this induces a real structure on $\mathcal{F}$. Letting $B_\R:=B\cap G_\R$, we get an identification of the real points $\mathcal{F}(\R)\simeq G_\R/B_\R$. For instance, if $K=SU(n)$, then $G=SL_n(\C)$ and $G_\R=SL_n(\R)$. We can take $B$ to be the subgroup of upper-triangular matrices in $SL_n(\C)$ and then the flag variety $\mathcal{F}_n:=SL_n(\C)/B$ is the classical variety of flags in $\C^n$. The Iwasawa decomposition in this case is the Gram--Schmidt process and we have diffeomorphisms
\[\mathcal{F}_n(\R)=SL_n(\R)/B_\R\simeq SO(n)/S(O(1)^n)\simeq O(n)/O(1)^n,\]
where $S(O(1)^n)$ (resp. $O(1)^n$) is the subgroup of diagonal matrices in $SO(n)$ (resp. in $O(n)$), which is isomorphic to $(\Z/2\Z)^{n-1}$ (resp. to $(\Z/2\Z)^n$).

Observe now that the compactness of the group $K$ ensures that it admits a bi-invariant Riemannian metric, yielding a Riemannian structure (called \textit{normal homogeneous}) on the flag manifold $\mathcal{F}$, as well as on $\mathcal{F}(\R)$, by restriction. Then, the action of $W$ on $\mathcal{F}$ (and on $\mathcal{F}(\R)$) is isometric. This fact is a first clue suggesting that the normal homogeneous Riemannian geometry on $\mathcal{F}$ may contain information on the combinatorics of the action of $W$ and, hopefully, could even provide a $W$-equivariant cell structure on $\mathcal{F}$.

In \cite{chirivi-garnier-spreafico}, inspired by methods for determining cellular structures on spherical space forms developed in a series of works initiated by \cite{mms13}, such a structure is constructed on the flag manifold $\mathcal{F}_3(\R)\simeq O(3)/O(1)^3$ of $SL_3(\R)$, equivariant with respect to the Weyl group $W=\Sym_3$. This is done by considering the free action of the \emph{binary octahedral group} $\mathcal{O}$ (of order 48) on the 3-sphere $\Sph^3$ and by noticing that this action gives the same information as the action of $\Sym_3$ on $\mathcal{F}_3(\R)$ (i.e. we have a diffeomorphism $\Sph^3/\mathcal{O}\simeq\mathcal{F}_3(\R)/\Sym_3$). The manifold $\mathcal{F}_3(\R)$ carries two natural metrics, which are proportional: the bi-invariant one inherited from $SO(3)$ and the one induced by modding out the (standard) round metric on $\Sph^3$ by the quaternion group $\mathcal{Q}_8$; see \S \ref{biinvariance} for more details. This allows us to interpret the cells of the $\Sym_3$-equivariant cell structure from \cite{chirivi-garnier-spreafico} as unions of open geodesics. For example, the 1-cells are (minimal) geodesic arcs between (the class of) $1$ and the reflections of $\Sym_3$, seen as points of $\mathcal{F}_3(\R)$; see Remark \ref{why}.

This suggests that an $\Sym_3$-equivariant cell decomposition can be obtained intrinsically, using only a normal homogeneous metric (a purely Lie-theoretical object) and not the special fact that $\mathcal{F}_3(\R)$ is a spherical space form. This will be the goal of the second part of this work.

In the more general case where $W\le\mathrm{Isom}(M)$ is a discrete isometry group of a connected complete Riemannian manifold $(M,g)$, we consider the \emph{Dirichlet--Voronoi domain}
\[\mathcal{DV}:=\{x\in M~;~\forall w\in W,~d(x_0,x)\le d(wx_0,x)\},\]
where $d$ is the geodesic distance on $M$ and $x_0\in M$ is a regular point. Such subsets were originally introduced by Dirichlet in the 1850s to study Fuchsian groups acting on the hyperbolic plane \cite{ratcliffe_hyperbolic}. More precisely, given a discrete subgroup $\Gamma\subset PSL_2(\R)$, the Dirichlet domain is a fundamental polyhedron for the $\Gamma$-action, whose translates tessellate the plane, equivariantly with respect to $\Gamma$. Moreover, such domains are algorithmically computable \cite{ratcliffe_hyperbolic,page_kleinian,voigt_fundomfuchs}, making them a useful explicit tool. 

This construction has since been extended to the Riemannian setting \cite{boissonnat-sharir-tagansky-yvinec,curved_voronoi_diagrams,boissonnat-rouxel-wintraecken}, with some surprising consequences: for example, one can find in \cite{chapron_phd} a probabilistic proof of the Gauss--Bonnet theorem, using Voronoi domains. Voronoi diagrams and related Delaunay triangulations also find applications in computer graphics and visualization, where various computation techniques were developed \cite{dyer-zhang-moller,NWTD}.

Regarding the construction of general Delaunay triangulations of Riemannian manifolds from Voronoi cells, first results were announced in \cite{leibon-letscher}, but there were later proven wrong \cite{boissonnat-dyer-ghosh_rep, boissonnat-dyer-ghosh-martynchuk}. As shown by the first general triangulation result \cite{boissonnat-dyer-ghosh17}, the situation is more intricate than expected and many strong hypotheses are required. The approach we follow in the present work is closer to the one of \cite{reflection_groups_on_manifolds}, where the case of a group acting by {\it reflections} on a Riemannian manifold is investigated.

The first part of this work is devoted to the exposition of general features of the set $\mathcal{DV}$, such as Proposition \ref{funddomdiscreteisometrygroup}, stating that $\mathcal{DV}$ is a path-connected fundamental domain for $W$ acting on $M$. Moreover, we detail a method for constructing equivariant cell structures from ordinary ones on $\mathcal{DV}$ that are admissible, in some sense. The precise statement is the following:
\begin{prop*}[\emph{Proposition \ref{cellsfromDV}}]
Suppose that $W$ acts freely on $M$ and that $\mathcal{DV}$ carries a regular CW structure such that, for any subset $I\subset W\setminus\{1\}$, the ``wall'' $Z_I:=\mathcal{DV}\cap\bigcap_{w\in I}w\mathcal{DV}$ is a (possibly empty) disjoint union of closed $|I|$-codimensional cells. Then, the $W$-translates of this structure form a $W$-equivariant CW structure on $M$.
\end{prop*}

Though rather elementary, this result will be key in the sequel, where we derive a new $\Sym_3$-equivariant CW structure on $\mathcal{F}_3(\R)$. Besides and as explained in Example \ref{lens_spaces}, even though the above result does not apply rigorously in this context, our method still gives insights on the construction of the classical CW structure of lens spaces \cite{hatcher}.

We next start the second part of this work by computing the injectivity radius of any (complex or real) flag manifold, equipped with a \textit{normal homogeneous metric}, i.e. one which is induced by a bi-invariant metric on the Lie group. Indeed, as stated in Proposition \ref{intDVcell}, the \textit{injectivity radius} of $(M,g)$ is an important invariant in our method. Such a metric on a flag manifold is unique up to a scalar and we choose to work with the metric induced by (minus) the Killing form on the Lie algebra and we call this metric \textit{standard homogeneous}. Specifically, after observing that the flag manifold $K/T=:\mathcal{F}_\Phi$ only depends on the root system of $K$, we prove the following result, which is essentially a consequence of Klingenberg's lemma \cite[Proposition 2.6.8]{klingenberg_riemann}:
\begin{theo*}[\emph{Theorem \ref{injrad}}]
The injectivity radius $\inj(\mathcal{F}_\Phi)$ of a standard homogeneous flag manifold $\mathcal{F}_\Phi$ is given by
\[\inj(\mathcal{F}_\Phi)=\pi\sqrt{h^\vee},\]
where $h^\vee$ is the dual Coxeter number of the (simple) root system $\Phi$. Moreover, this is the distance between 1 and any reflection of $W$, whose associated root is long.
\end{theo*}
In particular, we retrieve the B\"{o}ttcher--Wenzel inequality \cite{bottcher-wenzel} for skew-hermitian matrices and the Bloch--Iserles inequalities \cite{bloch-iserles} for skew-symmetric (real) matrices.

Back to the special case of the real flag manifold $\mathcal{F}_3(\R)=\mathcal{F}_{A_2}(\R)$, we prove that the ``radius'' of the Dirichlet--Voronoi domain $\mathcal{DV}_3\subset\mathcal{F}_3(\R)$ (i.e., the maximal distance from any of its points to its center) is smaller than the injectivity radius. Moreover, we show that the radius of $\mathcal{DV}_3$ is realized by exactly twenty-four points of $\partial\mathcal{DV}_3$; see Proposition \ref{DVisokforSO(3)} for more details. In particular, the interior $\interior{\mathcal{DV}}_3$ is a 3-cell and the twenty-four extremal points are some of the $0$-cells of the $\Sym_3$-equivariant cellular structure on $\mathcal{F}_3(\R)$, provided by Proposition \ref{cellsfromDV}. Moreover, the domain is combinatorially equivalent to a \emph{truncated cube}, yielding a (polyhedral) cellular decomposition of $\mathcal{DV}_3$; see Proposition \ref{ident_tcube} and Corollary \ref{ident_faces} for precise statements. We compute the cellular chain complex and obtain the following final result:
\begin{theo*}[\emph{Corollary \ref{indeedequivdec} and Theorem \ref{last_decomposition}}]
The Dirichlet--Voronoi domain $\mathcal{DV}$ is a fundamental domain for $\Sym_3$ acting on $\mathcal{F}_3(\R)$ and admits a cellular structure inducing an $\Sym_3$-equivariant cellular decomposition on $\mathcal{F}_3(\R)$, whose associated cellular homology chain complex is given by
\[\xymatrix{\Z[\Sym_3] \ar^{\partial_3}[r] & \Z[\Sym_3]^7 \ar^{\partial_2}[r] & \Z[\Sym_3]^{12} \ar^{\partial_1}[r] & \Z[\Sym_3]^6}\]
with boundaries
\[\partial_1=\left(\begin{smallmatrix}0 & 0 & 0 & 0 & 0 & s_\beta & -s_\beta & 0 & 0 & -1 & 0 & 1 \\ 0 & 0 & 0 & 1 & -1 & 0 & 0 & 0 & s_{\beta}s_\alpha & 0 & 0 & -1 \\ -w_0 & 0 & 0 & 0 & 0 & 0 & 0 & w_0 & 0 & s_\beta & -w_0 & 0 \\ s_\beta s_\alpha & -s_\beta s_\alpha & 0 & 0 & s_\alpha & -s_\alpha & 0 & 0 & 0 & 0 & 0 & 0 \\ 0 & s_\beta s_\alpha & -s_\beta s_\alpha & 0 & 0 & 0 & 0 & 0 & -w_0 & 0 & w_0 & 0 \\ 0 & 0 & 1 & -1 & 0 & 0 & s_\beta & -s_\beta & 0 & 0 & 0 & 0\end{smallmatrix}\right),\]

\[\partial_2=\left(\begin{smallmatrix}
1 & 0 & w_0 & 0 & 0 & 0 & -w_0 \\
1 & -s_\alpha s_\beta & 0 & 0 & -1 & 0 & 0 \\
1 & s_\beta & 0 & -s_\beta & 0 & 0 & 0 \\
1 & 0 & s_\alpha & 0 & 0 & -1 & 0 \\
1 & 0 & -1 & 0 & -w_0 & 0 & 0 \\
1 & s_\alpha & 0 & 0 & 0 & 0 & -1 \\
1 & -1 & 0 & 0 & 0 & -s_\beta & 0 \\
1 & 0 & -s_\beta s_\alpha & -1 & 0 & 0 & 0 \\
0 & -1 & -w_0 & 0 & -s_\beta & 0 & 0 \\
0 & s_\beta & 1 & 0 & 0 & 0 & -s_\beta \\
0 & w_0 & -w_0 & -1 & 0 & 0 & 0 \\
0 & -s_\beta s_\alpha & 1 & 0 & 0 & -1 & 0\end{smallmatrix}\right),~~
\partial_3:=\left(\begin{smallmatrix}
1-s_\alpha \\
1-s_\beta \\
1-w_0 \\
1-s_\beta s_\alpha \\
1-s_\alpha s_\beta \\
1-s_\beta s_\alpha \\
1-s_\alpha s_\beta\end{smallmatrix}\right),\]
where $s_\alpha$ and $s_\beta$ are the simple reflections of $\Sym_3$ and $w_0:=s_\alpha s_\beta s_\alpha$ is its longest element.
\end{theo*}

In the Appendix \ref{truncated_cube} we provide some figures describing the combinatorics of the polytopal structure of the domain $\mathcal{DV}_3\subset\mathcal{F}_3(\R)$.

In general, we make the following conjecture, roughly saying that Proposition \ref{cellsfromDV} applies to all real flag manifolds:
\begin{conjecture*}
{We endow $\mathcal{F}_\Phi=K/T$ with any normal homogeneous metric.}
\begin{enumerate}
\item The Dirichlet--Voronoi domain $\mathcal{DV}$ associated to $\mathcal{F}_\Phi$ and $W$ is included in the open (geodesic) ball centered at $1$ and of radius $\inj(\mathcal{F}_\Phi)$.
\item If $1\notin I\subset W$ is a subset, then the wall $\mathcal{F}_\Phi(\R)\cap Z_I$ is a (possibly empty) union of $(N-|I|)$-cells, where $N:=\dim\mathcal{F}_\Phi(\R)=\tfrac12(\dim K-\rk K)$.
\end{enumerate}
\end{conjecture*}
These are to be addressed in future works.

\part{Generalities on Dirichlet--Voronoi domains in Riemannian manifolds and related cell structures}

In this first part, we outline a general method for constructing cell structures on Riemannian manifolds, equivariant with respect to a finite group of isometries, using Dirichlet--Voronoi domains, in view of applying it to real flag manifolds. 

The first section gives general results on Dirichlet--Voronoi domains. In particular, we exhibit a condition on the \textit{injectivity radius} of the manifold under which the interior of such a domain is a cell; see Propositions \ref{funddomdiscreteisometrygroup} and \ref{intDVcell}. We also give a useful criterion (Lemma \ref{isolated_max}) for a geodesic ball to contain the domain.

Recall that in the case of Euclidean, spherical and hyperbolic spaces \cite[Theorem 6.8.1]{ratcliffe_hyperbolic}, if the Dirichlet domain $\mathcal{P}$ is a (convex) fundamental polyhedron for the action of a discrete group $\Gamma$, then $\{g\mathcal{P},~g\in\Gamma\}$ is a tessellation of the ambient space and, in particular, the ``walls'' $\mathcal{P}\cap g\mathcal{P}$ tessellate the boundary $\partial\mathcal{P}$ of $\mathcal{P}$. In our setting, once we know that the open Dirichlet domain is a top-cell, we look for a cell structure on its boundary, with cells given by the (connected components of the) walls. In Proposition \ref{cellsfromDV}, we prove that some compatibility assumptions on a (regular) CW structure on the domain allow us to induce an equivariant cell structure on the whole manifold. Of course, such a structure cannot be expected to exist in general, as Example \ref{lens_spaces} of lens spaces shows\footnote{However, a cell structure on all lens spaces can still be ``obtained'' in this way; see Example \ref{lens_spaces}.}. Nevertheless, using the quaternion algebra, in \S \ref{O(3)/O(1)3} we will prove that this is indeed the case for the flag manifold $\mathcal{F}_3(\R)$ of $SL_3(\R)$, thus leading to a new $\Sym_3$-equivariant cell structure on it.

\section{Definitions and general properties}

In this section, we fix $(M,g)$ a connected complete Riemannian manifold, with geodesic distance $d$, and $W\le\mathrm{Isom}(M)$ a discrete subgroup of the isometry group of $(M,g)$. By a classical result (see \cite[Lemma 2.1]{reflection_groups_on_manifolds} for instance), this means that each $W$-orbit is discrete as a subset of $M$. We fix also $x_0\in M$ a \emph{regular point}, i.e. a point with trivial stabilizer.

Inspired by the study of Fuchsian groups and following the approach of \cite{reflection_groups_on_manifolds}, we consider the Dirichlet--Voronoi domain of $W$ acting on $\mathcal{F}$.

\begin{definition}\label{dirichlet}
Let $x_0\in M$ be a regular point.
\begin{enumerate}[label=$\bullet$]
\item The \emph{Dirichlet--Voronoi domain} centered at $x_0$ is the following subset of $M$:
\[\mathcal{DV}:=\{x\in M~;~\forall w \in W,~d(x_0,x)\le d(wx_0,x)\}.\]
\item For $w\in W$, we denote by $H_w$ the \emph{dissecting hypersurface}
\[H_w:=\{x\in M~;~d(x_0,x)=d(wx_0,x)\}\]
and by $Z_w$ the (maximal) \emph{wall}
\[Z_w:=\mathcal{DV}\cap H_w=\mathcal{DV}\cap w\mathcal{DV}.\]
\item For $I\subset W$, the \emph{$I$-wall} $Z_I$ is the (possibly empty) intersection of the maximal walls in $I$, i.e.
\[Z_I:=\bigcap_{w\in I}Z_w=\mathcal{DV}\cap\bigcap_{w\in I}w\mathcal{DV}.\]
\end{enumerate}
\end{definition}

If the action is free, we may consider the orbit space $M/W$, equipped with the quotient metric $g/W$ and geodesic distance $d_{M/W}$. Then, the set $\mathcal{DV}$ can be interpreted as the set of elements $x\in M$ realizing the distance of their orbit: $d(x_0,x)=d_{M/W}(Wx_0,Wx)$.

\begin{rem}\label{care_with_convexity}
As already mentioned, the domains $\mathcal{DV}$ as defined above are mainly studied for hyperbolic manifolds (see \cite{bowditch_hyperbolic}) or more generally for manifolds with constant sectional curvature (see \cite[\S 6.6]{ratcliffe_hyperbolic}). This is because we want $\mathcal{DV}$ to be a \emph{fundamental polyhedron} for $W$ acting on $M$ and in particular, geodesically convex. In the case of flag manifolds, the curvature is no longer constant and one has to be careful with the meaning of convexity, because minimal geodesics are not unique in general. As we shall see, a relevant notion to introduce regarding this matter is the \emph{injectivity radius} $\inj_{x_0}(M)$ of $M$ at $x_0$ (see \cite[\S 6.2]{lee_riema} or \cite[Definition 2.116]{gallot-hulin-lafontaine}).
\end{rem}

It follows immediately from the above definition that for $w\in W$, the subset $M\setminus H_w$ is the disjoint union of the two open subsets $\{x\in M~;~d(x_0,x)<d(wx_0,x)\}$ and $\{x\in M~;~d(x_0,x)>d(wx_0,x)\}$ and moreover, the interior of $\mathcal{DV}$ is the connected component of $M\setminus \bigcup_{1\ne w\in W}H_w$ containing $x_0$.

It is reasonable to expect $\mathcal{DV}$ to be a fundamental domain for $W$ acting on $M$. This is indeed the case, but we shall need a technical preliminary result on the behavior of the hypersurfaces $H_w$ with respect to minimal geodesics. This result appears in \cite{reflection_groups_on_manifolds}, we reproduce it here for the sake of self-containment. One should be careful with the terminology: though we call the $H_w$'s ``hypersurfaces'', they are not necessarily submanifolds of $M$.

\begin{lem}[\emph{\cite[Lemma 2.2]{reflection_groups_on_manifolds}}]\label{interonlyonce}
Let $y_0,y_1\in M$ be distinct points of $M$ and consider the hypersurface $H:=H_{y_0,y_1}=\{x\in M~;~d(x,y_0)=d(x,y_1)\}$. If $x\in H$, then every minimal geodesic from $y_0$ to $x$ meets $H$ only at $x$.
\end{lem}

\begin{proof}
Let $\gamma_0$ be a minimal geodesic parametrized by arc-length such that $\gamma_0(0)=y_0$ and $\gamma_0(\ell)=x$, where $\ell=d(y_0,x)$ and suppose for contradiction that $\gamma_0(t)\in H$ for some $t<\ell$. We compute
\[d(x,y_1)=d(x,y_0)=d(y_0,\gamma_0(t))+d(\gamma_0(t),x)=d(y_1,\gamma_0(t))+d(\gamma_0(t),x),\]
and so we are in the case of equality in the triangular inequality. Let $\gamma_1$ be a minimal geodesic from $\gamma_0(t)$ to $y_1$ and let $\widetilde{\gamma_0}$ be the curve $s\mapsto \gamma_0(t+s)$ for $0\le s\le \ell-t$. The situation can be visualized in Figure \ref{illus_iteronlyonce}.

Then, the concatenation $\gamma_2:=\gamma_1\ast \widetilde{\gamma_0}^{-1}$ is a piecewise smooth curve from $x$ to $y_1$ satisfying
\[\leng(\gamma_2)=\leng(\gamma_1)+\leng(\widetilde{\gamma_0})=\ell=d(x,y_1),\]
where $\leng$ denotes the length of a curve. By \cite[Chapter 3, Corollary 3.9]{do_carmo_riemannian}, this implies that $\gamma_2$ is in fact a (smooth) minimal geodesic from $x$ to $y_1$, which coincides with the geodesic $\widetilde{\gamma_0}^{-1}$ on a non-empty interval. By the Picard-Lindel\"{o}f theorem, this implies $\gamma_2(s)=\gamma_0(\ell-s)$ for $0\le s \le \ell$ and thus $y_0=\gamma_0(0)=\gamma_2(\ell)=y_1$, a contradiction.
\end{proof}
\begin{center}
\begin{figure}[h!]
\centering
\begin{tikzpicture}[scale=1.8]
	\coordinate (yz) at (-1,-1);
	\coordinate (yu) at (1,1);
	\coordinate (bha) at (-0.2,1.3);
	\coordinate (bhb) at (1,-0.5);
	\coordinate (z) at (-0.35,0.5);
	\coordinate (x) at (0.3,-0.3);
	
	\draw (yz) node[below]{$y_0$};
	\draw (yu) node[above]{$y_1$};
	\draw (x) node[below left]{$x$};
	\draw (z) node[left]{$\gamma_0(t)$};
	\draw (0.35,0.4) node{$\widetilde{\gamma_0}$};
	\draw (bhb) node[right]{$H$};
	\fill[fill=black] (yz) circle (1pt);
	\fill[fill=black] (yu) circle (1pt);
	\fill[fill=black] (x) circle (1pt);
	\fill[fill=black] (z) circle (1pt);
	
	\draw[xshift=4cm,dashed,thick] plot [smooth,tension=1] coordinates {(bha) (z) (x) (bhb)};
	\draw[thick] plot [smooth,tension=2] coordinates {(yz) (z) (x)} [arrow inside={end=to,opt={black,scale=1}}{0.25,0.8}];
	\draw[thick] plot [smooth] coordinates {(z) (0.2,1) (yu)} [arrow inside={end=to,opt={black,scale=1}}{0.25,0.75}];
	
	\draw (-0.85,-0.15) node[left]{$\gamma_0$};
	\draw (0.2,1) node[above]{$\gamma_1$};
\end{tikzpicture}
\caption{Illustration of the proof of Lemma \ref{interonlyonce}.}\label{illus_iteronlyonce}
\end{figure}
\end{center}

Another interesting feature of $\mathcal{DV}$ is that it is path-connected. More precisely, we have the following result:

\begin{lem}\label{DVisstarshaped}
The domain $\mathcal{DV}$ is \emph{geodesically star-shaped} with respect to $x_0$, meaning that for every $x\in \mathcal{DV}$ and any minimal geodesic $\gamma : [0,1]\to M$ from $x_0$ to $x$, we have $\gamma(t)\in\mathcal{DV}$ for every $t\in[0,1]$. In particular, $\mathcal{DV}$ is path-connected.
\end{lem}

\begin{proof}
Let $t\in[0,1]$ and $w\in W$. We write
\[\begin{array}[t]{ll}
d(x_0,w\gamma(t)) &\begin{array}[t]{lr}
					=d(x_0,w\gamma(t))+d(w\gamma(t),wx)-d(\gamma(t),x) & \text{($w$ is an isometry)} \\
					\ge d(x_0,wx)-d(\gamma(t),x) & \text{(triangular inequality)}\\
					\ge d(x_0,x)-d(\gamma(t),x) & \text{($x_0\in\mathcal{DV}$)}\\
					=d(x_0,\gamma(t)) & \text{($\gamma$ is minimal)}\end{array}\\
\end{array}\]
and therefore we have $\gamma(t)\in\mathcal{DV}$, as required.
\end{proof}

\begin{prop}\label{funddomdiscreteisometrygroup}
The Dirichlet--Voronoi domain $\mathcal{DV}$ is a geodesically star-shaped fundamental domain for $W$ acting on $M$.
\end{prop}

\begin{proof}
Obviously we have $M=\bigcup_{w\in W}w\mathcal{DV}$. On the other hand, if $x\in \mathcal{DV}\cap w\mathcal{DV}=Z_w$ for some $w\in W\setminus\{1\}$ and if $B=B(x,\delta)$ is a small (geodesic) ball centered at $x$ with radius $\delta>0$ included in $Z_w$, then $B\subset H_w$. However, if we denote by $\gamma$ a minimal geodesic from $x_0$ to $x$ (parametrized by arc-length) and if $\ell:=d(x_0,x)=\leng(\gamma)$, then $d(\gamma(t),x)=\ell-t<\delta$ for $t>\ell-\delta$ and thus $\gamma(t)\in H_w$ for $\ell-\delta<t\le\ell$, contradicting Lemma \ref{interonlyonce}. Therefore, $Z_w$ has empty interior.
\end{proof}

\section{Injectivity questions}

We intend to use the domain $\mathcal{DV}$ to build a $W$-equivariant CW structure on $M$. Of course, this is too much to ask in the general setting, as the walls $Z_w$ are not necessarily cells. For example, let the cyclic group $C_2=\{1,s\}$ act on $\Sph^2$ via the antipode. Then the wall $Z_s$ is a circle. We can however take a new Dirichlet--Voronoi domain for $C_2$ acting on $Z_s$, leading to an actual $C_2$-CW structure on $\Sph^2$; see Example \ref{examSU2}. This will hopefully lead to a general construction of $W$-equivariant cell structures on complex flag manifolds. In the present work though, we stick to the case of real flag manifolds.

The first feature to ask for is that the interior of $\mathcal{DV}$ should itself be a cell. To ensure this, we have to control the size of $\mathcal{DV}$. The following elementary result, which doesn't seem to appear in the literature (as far as the author knows), can be compared with \cite[Lemma B.4]{boissonnat-dyer-ghosh_rep} and \cite[Theorem 1]{dyer-vegter-wintraecken}.

Before stating the result, we introduce some notation: for any $x\in B(x_0,\inj_{x_0}(M))$ we denote by $\gamma_x$ the \emph{unique} minimal geodesic from $x_0$ to $x$, extended to $\R$ by completeness of $M$. The geodesic $\gamma_x$ is defined by $\gamma_x(s)=\mathrm{Exp}_{x_0}\!\left({su}/{\|u\|}\right)$ for any $s\in\R$, where $u=\dot{\gamma_x}(0)\in T_{x_0}M$ is such that $\Exp_{x_0}(u)=x$; in particular, $\|u\|=d(x_0,x)$.

\begin{prop}\label{intDVcell}
If $0<\rho\le\inj_{x_0}(M)$ is such that the open fundamental domain $\stackrel{\circ}{\mathcal{DV}}$ is included in the geodesic ball $B(x_0,\rho)$, then $\stackrel{\circ}{\mathcal{DV}}$ is a $\dim(M)$-cell.
\end{prop}

\begin{proof}
Since $\stackrel{\circ}{\mathcal{DV}}\subset B(x_0,\rho)$, we have $\mathcal{DV}\subset\overline{B(x_0,\rho)}$. Thus, by Lemma \ref{interonlyonce}, for each $y\in B_{T_{x_0}M}(0,\rho)$ there is a unique $0<\rho_y\le\rho$ such that $\gamma_{\Exp(y)}(\rho_y)=\Exp_{x_0}(\rho_yy/\|y\|)$ is in $\partial\mathcal{DV}$. Therefore, the element $\Exp_{x_0}(\rho_yy/\rho)$ is in the interior of $\mathcal{DV}$, so that the assignment
\[\begin{array}{ccc}
B_{T_{x_0}M}(0,\rho) & \stackrel{\tiny{\Psi}}\longto & \stackrel{\circ}{\mathcal{DV}} \\ y & \longmapsto & \Exp_{x_0}(\rho_yy/\rho)\end{array}\]
defines a continuous map. Conversely, if $x\in~\stackrel{\circ}{\mathcal{DV}}$, then there is a unique $0<\rho_x\le\rho$ such that $\gamma_x(\rho_x)\in\partial\mathcal{DV}$. Letting $\ell_x:=d(x_0,x)<\rho_x$, we obtain a continuous map
\[\begin{array}{ccc}
\stackrel{\circ}{\mathcal{DV}} & \longto & B_{T_{x_0}M}(0,\rho) \\ x & \longmapsto & \mathrm{Exp}_{x_0}^{-1}(\gamma_x(\rho\ell_x/\rho_x))=\rho \Exp_{x_0}^{-1}(x)/\rho_x\end{array}\]
which is easily seen to be an inverse for $\Psi$.
\end{proof}

Another observation towards the construction of cell structures is the following: assume we can find some $0<\delta<\inj_{x_0}(M)$ such that $\mathcal{DV}\subset \overline{B(x_0,\delta)}$. Then, the exponential map $\mathrm{Exp}={\rm Exp}_{x_0}$ realizes a homeomorphism $\overline{B(x_0,\delta)}\simeq\overline{B_{T_{x_0}M}(0,\delta)}=\mathbb{B}^{N}$, where $\mathbb{B}^{N}$ is the Euclidean $N$-ball, with $N:=\dim M$. Thus, we can project the ``cells'' onto the bounding sphere $\partial\overline{B(x_0,\delta)}=S(x_0,\delta)\simeq S_{T_{x_0}M}(0,\delta)$, just as in the case of binary polyhedral groups (see \cite{chirivi-spreafico} and \cite{chirivi-garnier-spreafico}). More precisely, we consider the map
\[\begin{array}{ccccc}
\pi_\delta & : & \overline{B(x_0,\delta)}\setminus\{x_0\} & \longtwoheadrightarrow & S_{T_{x_0}M}(0,\delta)\simeq\Sph^{N-1} \\ & & x & \longmapsto & \delta\frac{\mathrm{Exp}^{-1}(x)}{\|\mathrm{Exp}^{-1}(x)\|}\end{array}\]
Geometrically, this can be interpreted in terms of geodesics: for $x\in\overline{B(x_0,\delta)}$, as $d(x_0,x)<\inj_{x_0}(M)$, there is a unique minimal geodesic $\gamma$ from $x_0$ to $x$, which we extend until it meets the sphere $S(x_0,\delta)$; the preimage of this intersection point under $\mathrm{Exp}$ is $\pi_\delta(x)$. Now, if we have a wall $Z_w$ and if $y=\pi_\delta(x)$ for $x\in Z_w$, then the unique minimal geodesic from $x_0$ to $\mathrm{Exp}(y)$ intersects $Z_w$ in at least one point, hence in a unique point thanks to Lemma \ref{interonlyonce}. Therefore, $x=\pi_\delta^{-1}(y)$ is well-defined and thus $\pi_\delta$ restricts to a homeomorphism $Z_w \stackrel{\tiny{\sim}}\to \pi_\delta(Z_w)\subset\Sph^{N-1}$. We may glue these homeomorphisms together to obtain a global homeomorphism
\[\pi_\delta : \partial\mathcal{DV}\stackrel{\tiny{\sim}}\longto \Sph^{N-1}.\]
One should notice that the walls $Z_w$ are not connected in general: they are rather expected to be finite unions of cells. Though this fails in general --- see Examples \ref{lens_spaces} for lens spaces and \ref{examSU2} for complex flag manifolds --- we shall see that this is the case for the flag manifold of $SL_3(\R)$. 

We finish this section with a technical lemma that helps us find a bound on $\delta>0$ such that $\mathcal{DV}\subset \overline{B(x_0,\delta)}$, assuming (a strong form of) the injectivity radius condition and that the acting group $W$ is finite:

\begin{lem}[``No antenna lemma'']\label{isolated_max}
Let $W\le\mathrm{Isom}(M)$ be finite, with associated Dirichlet--Voronoi domain $\mathcal{DV}$ and assume that $\mathcal{DV}\subset B(x_0,\rho)$ for some $0<\rho<\inj_{x_0}(M)$. If $0<\delta<\rho$ is such that the intersection $\mathcal{DV}\cap S(x_0,\delta)$ of $\mathcal{DV}$ with the sphere of radius $\delta$ consists only of isolated points, then $\mathcal{DV}\subset\overline{B(x_0,\delta)}$.
\end{lem}

\begin{proof}
If there is some $z\in\mathcal{DV}$ such that $d(x_0,z)>\delta$, then $x:=\gamma_z(\delta)\in\mathcal{DV}\cap S(x_0,\delta)$ and thus for any $0<\varepsilon<d(x,z)$, the element $\gamma_z(\delta+\varepsilon/2)$ is in $\mathcal{DV}\cap{}^{\complement}{B(x_0,\delta)}\cap B(x,\varepsilon)$. We will prove however that this set is empty for $\varepsilon>0$ sufficiently small.

As $\mathcal{DV}\cap S(x_0,\delta)$ consists of isolated points, we may choose $0<\varepsilon<\tfrac{\rho-\delta}{2}$ such that 
\[\mathcal{DV}\cap S(x_0,\delta)\cap B(x,2\varepsilon)=\{x\}.\]
Suppose for contradiction that $y\in\mathcal{DV}\cap{}^{\complement}{B(x_0,\delta)}\cap B(x,\varepsilon)$. We define $d_y:=d(x_0,y)\ge\delta$ and compute
\[d(x,\gamma_y(\delta))\le d(x,y)+d(y,\gamma_y(\delta))<\varepsilon+d(\gamma_y(d_y),\gamma_y(\delta))=\varepsilon+d_y-\delta\le\varepsilon+d(x,y)<2\varepsilon,\]
so $\gamma_y(\delta)\in\mathcal{DV}\cap S(x_0,\delta)\cap B(x,2\varepsilon)$ and $\gamma_y(\delta)=x=\gamma_z(\delta)$. Therefore, $\gamma_y=\gamma_z$ as there is only one minimal geodesic from $x_0$ to $x$, since $x\in\mathcal{DV}\subset B(x_0,\inj_{x_0}\!M)$. This proves that
\[\mathcal{DV}\cap{}^{\complement}{B(x_0,\delta)}\cap B(x,\varepsilon)=\gamma_z([\delta,\delta+\varepsilon[).\]
The situation (which we want to prove impossible) is depicted in Figure \ref{no_antenna}, giving the lemma its name.
\begin{figure}
\centering
\begin{tikzpicture}[scale=2.5]
	\coordinate (z) at (0,0);
	\coordinate (a) at (0,1);
	\coordinate (b) at (-0.809,0.588);
	\coordinate (c) at (1,0);
	\coordinate (d) at (0.7071,-0.7071);
	\coordinate (e) at (-0.39,-0.92);
	\coordinate (Z) at (0.1,1.3);
	\coordinate (Zp) at (0.17,1.4);
	
	\draw (z) circle (1);
	\fill[fill=black] (a) circle (0.75pt);
	\fill[fill=black] (z) circle (0.75pt);
	\fill[fill=black] (b) circle (0.75pt);
	\fill[fill=black] (c) circle (0.75pt);
	\fill[fill=black] (d) circle (0.75pt);
	\fill[fill=black] (e) circle (0.75pt);
	\fill[fill=black] (Zp) circle (0.75pt);
	\draw (a) circle (0.315);
	
	\draw[dashed,thick] plot [smooth,tension=1] coordinates {(e) (-0.3,-0.08) (b)};
	\draw[dashed,thick] plot [smooth,tension=1] coordinates {(e) (0.127,-0.651) (d)};
	\draw[dashed,thick] plot [smooth,tension=1] coordinates {(d) (0.6,-0.25) (c)};
	\draw[dashed,thick] plot [smooth,tension=1] coordinates {(c) (0.4,0.4) (a)};
	\draw[dashed,thick] plot [smooth,tension=1] coordinates {(a) (-0.303,0.6) (b)};
	
	\draw (a) node[above left]{$x$};
	\draw (z) node[below]{$x_0$};
	\draw (Zp) node[right]{$z$};
	\draw (-0.1,0.6) node[below]{$\gamma_z$};
	\draw (-0.15,1.3) node[left]{$S(x,\varepsilon)$};
	\draw (0.66,0.66) node[above right]{$S(x_0,\delta)$};
	\draw (-0.22,0) node[below left]{$H_{w_1}$};
	\draw (0.13,-0.651) node[below]{$H_{w_2}$};
	\draw (0.55,-0.16) node[below right]{$H_{w_3}$};
	\draw (0.4,0) node{$\mathcal{DV}$};
	
	\draw[thick] plot [smooth,tension=1] coordinates {(z) (0.1,0.3) (-0.04,0.6) (a) (Z)} [arrow inside={end=to,opt={black,scale=1}}{0.35,0.9}];;
	\draw[ultra thick] plot [smooth,tension=1] coordinates {(a) (0.038,1.15) (Z)};
	\draw[thick] plot [smooth,tension=1] coordinates {(a) (0.038,1.15) (Z) (Zp)};
\end{tikzpicture}
\caption{The thick curve represents the forbidden antenna $\gamma_z([\delta,\delta+\varepsilon[)$.}
\label{no_antenna}
\end{figure}

On the other hand, by Lemma \ref{interonlyonce}, if we have $d(x_0,\gamma_z(t))=d(wx_0,\gamma_z(t))$ for some $1\ne w\in W$ and some $\delta<t<\delta+\varepsilon$, then $t=:t_w$ is unique (we set $t_w:=\delta+\varepsilon/2$ in other cases) and therefore, if $\delta<s<t_w$, then $d(x_0,\gamma_z(s))<d(wx_0,\gamma_z(s))$. Since $W$ is finite, we may choose $t_0$ such that $\delta<t_0<\min_{1\ne w\in W}t_w<\delta+\varepsilon$ and we then have $\gamma_z(t_0)\in\interior{\mathcal{DV}}$.

For a unit vector $v\in T_{x_0}M$, we let $\gamma^v : \R\to M$ be the geodesic $s\mapsto \mathrm{Exp}_{x_0}(sv)$. Defining $S_T(0,1):=S_{T_{x_0}M}(0,1)$, the set
\[\{v\in S_T(0,1)~;~\gamma^v(t_0)\in\interior{\mathcal{DV}}\}=\{v\in S_T(0,1)~;~d(x_0,\gamma^v(t_0))<d(wx_0,\gamma^v(t_0)),~\forall w\in W\setminus\{1\}\}\]
is an open neighbourhood of $u:=\dot{\gamma_z}(0)$ in $S_T(0,1)$. Hence, we may choose $0<\eta<1$ such that
\[\forall v\in S_T(0,1),~\|u-v\|<\eta~\Longrightarrow~\gamma^v(t_0)\in\interior{\mathcal{DV}}.\]
Since the map $\mathrm{Exp}_{x_0}$ is continuous, by shrinking $\eta$ if needed and since $\gamma_z(t_0)\in B(x,\varepsilon)$, we may assume that $\gamma^v(t_0)\in B(x,\varepsilon)$ for $\|u-v\|<\eta$. As $d(x_0,\gamma^v(t_0))=t_0>\delta$, we obtain
\[\forall v\in S_T(0,1),~\|u-v\|<\eta~\Longrightarrow~\gamma^v(t_0)\in \interior{\mathcal{DV}}\cap {}^\complement{B(x_0,\delta)}\cap B(x,\varepsilon)\subset\gamma_z([\delta,\delta+\varepsilon[).\]

In particular, we can find some $v\ne\pm u$ such that $\mathrm{Exp}_{x_0}(t_0v)\stackrel{\tiny{\text{df}}}=\gamma^v(t_0)=\gamma_z(s)\stackrel{\tiny{\text{df}}}=\Exp_{x_0}(su)$ for some $\delta\le s<\delta+\varepsilon$. But since $t_0$ and $s$ are no greater than $\delta+\varepsilon<\tfrac{1}{2}(\rho+\delta)<\rho<\inj_{x_0}\!(M)$, this implies that $t_0v=su$. This would mean that $u$ and $v$ are collinear and on the same sphere, a contradiction.
\end{proof}

\begin{rem}
In the case where $\inj_{x_0}\!(M)<\infty$, the existence of some number $0<\rho<\inj_{x_0}\!(M)$ such that $\mathcal{DV}\subset B(x_0,\rho)$ is equivalent to the assumption that $\mathcal{DV}\subset B(x_0,\inj_{x_0}\!(M))$.
\end{rem}

\section{Equivariant cell structures from compatible structures on $\mathcal{DV}$}

In this section, we give some preliminary results on how to build a $W$-equivariant cell structure on $M$ from a compatible cell structure on $\mathcal{DV}$. We need a preliminary result, giving some information on the partial action of $W$ on walls. The reader may observe that the first part of the following lemma appears in \cite[Corollary 1.4]{voigt_fundomfuchs}, where it is established for Fuchsian groups:
\begin{lem}\label{partactonwalls}
For $v,w\in W$, we have 
\[w^{-1}Z_w=Z_{w^{-1}},\]
as well as the \emph{cocycle inclusion}\footnote{This terminology is inspired by group cohomology: recall that for a finite group $G$ and a $\Z[G]$-module $M$, a 1-cocycle is a map $f : G\to M$ such that $f(gh)=f(g)+gf(h)$ for $g,h\in G$.}
\[Z_v\cap vZ_w\subset Z_{vw}.\]
\end{lem}
\begin{proof}
We first prove the equality. Take $z\in Z_w=\mathcal{DV}\cap H_w$. We calculate
\[d(w^{-1}x_0,w^{-1}z)=d(x_0,z)=d(wx_0,z)=d(x_0,w^{-1}z),\]
proving that $w^{-1}z\in H_{w^{-1}}$. To show that $w^{-1}z\in\mathcal{DV}$, take $u\in W$ and write
\[d(ux_0,w^{-1}z)=d(wux_0,z)\ge d(x_0,z)=d(wx_0,z)=d(x_0,w^{-1}z).\]
We thus have $w^{-1}Z_w\subseteq Z_{w^{-1}}$ and replacing $w$ by $w^{-1}$ yields the reverse inclusion. To prove the cocycle inclusion, if $x=vy\in H_v\cap vH_w$ with $y\in H_w$, then we have
\[d(vwx_0,x)=d(wx_0,v^{-1}x)=d(wx_0,y)=d(x_0,y)=d(vx_0,vy)=d(vx_0,x)=d(x_0,x)\]
and thus $x\in H_{vw}$, as required.
\end{proof}

\begin{prop}\label{cellsfromDV}
Suppose that $W$ acts freely on $M$ and that $\mathcal{DV}$ admits a regular CW structure
\[\mathcal{DV}=\coprod_{e\in\mathcal{E}}e\]
such that for a subset $1\notin I\subset W$, the $I$-wall $Z_I$ is a (possibly empty) disjoint union of closed cells of codimension $|I|$. If $e$ is an open cell in $Z_I$ and if $w\in W\setminus\{1\}$ is such that $w(e)\cap\mathcal{DV}\ne\emptyset$, then $w^{-1}\in I$ and $w(e)$ is an open cell in $wZ_I=Z_{w^{-1}}\cap Z_{wI\setminus\{1\}}$.

Moreover, the decomposition
\[M=\coprod_{\substack{e\in\mathcal{E} \\ w\in W}}we\]
is a $W$-equivariant CW structure on $M$.
\end{prop}
\begin{proof}
We start with the first statement: let $x\in we\cap\mathcal{DV}\ne\emptyset$. We have $w^{-1}x\in\mathcal{DV}\cap w^{-1}\mathcal{DV}=Z_{w^{-1}}$ and so 
\[w^{-1}x\in e\cap Z_{w^{-1}}\subset Z_{w^{-1}}\cap Z_I.\]
But if $w^{-1}\notin I$, this last intersection is at most a union of closed $(\dim(M)-|I|-1)$-cells, disjoint from $e$. Therefore, $w^{-1}\in I$ and using Lemma \ref{partactonwalls}, we get
\begin{align*}
wZ_I&=w\left(Z_{w^{-1}}\cap\bigcap_{w^{-1}\ne v\in I}Z_v\right)=wZ_{w^{-1}}\cap\bigcap_{v\ne w^{-1}}(wZ_{w^{-1}}\cap wZ_v)=Z_w\cap\bigcap_{v\ne w^{-1}}(Z_w\cap wZ_v) \\
&\subset Z_w\cap\bigcap_{v\ne w^{-1}}Z_{wv}=Z_w\cap Z_{wI\setminus\{1\}}.
\end{align*}
Conversely, using Lemma \ref{partactonwalls} again yields
\[w^{-1}(Z_w\cap Z_{wI\setminus\{1\}})=w^{-1}\left(Z_w\cap\bigcap_{w^{-1}\ne v\in I}Z_{wv}\right)=Z_{w^{-1}}\cap\bigcap_{v\ne w^{-1}}(Z_{w^{-1}}\cap w^{-1}Z_{wv})\subset Z_I\]
and we obtain
\[wZ_I=Z_w\cap Z_{wI\setminus\{1\}}.\]
As $e\subset Z_I$, we have $we\in Z_w\cap Z_{wI\setminus\{1\}}\subset\mathcal{DV}$ and since $e$ is the interior (in the $\dim(e)$-skeleton of $\mathcal{DV}$) of a connected component of $Z_I$, the translate $we$ is also the interior of a connected component of $Z_w\cap Z_{wI\setminus\{1\}}$, hence an open cell.

Let $e,e'\in\mathcal{E}$ be two cells and $w,w'\in W$ such that $we\cap w'e'\ne\emptyset$. We have to prove that $we=w'e'$. We have $e\cap (w^{-1}w'e')=w^{-1}(we\cap w'e')\ne\emptyset$ and thus $w^{-1}w'e'\cap\mathcal{DV}\ne\emptyset$. By the first statement, this implies that $w^{-1}w'e'\subset\mathcal{DV}$, so $e,w^{-1}w'e'\in\mathcal{E}$ and since $\mathcal{DV}=\coprod_ee$ is a cell decomposition, this yields $e=w^{-1}w'e'$, as required. Since $M=\bigcup_{w\in W}w\mathcal{DV}=\bigcup_{w,e}we$, this ensures that $M=\coprod_{w,e}we$ is indeed a CW structure. Moreover, it is clear that $W$ permutes the cells of $M$, so it remains to show that if $we=e$ for some $e\in\mathcal{E}$ and $w\in W$, then $w$ is the identity on $e$. As the action is free, this amounts to say that if $we=e$, then $w=1$. Choose a characteristic map $\interior{\psi} : \interior{\mathbb{B}}^k\stackrel{\tiny{\sim}}\longto e$. Since the CW structure on $\mathcal{DV}$ is regular, the map $\interior{\psi}$ extends to a homeomorphism $\psi : \mathbb{B}^k\stackrel{\tiny{\sim}}\longto\overline{e}$ and the map $\psi^{-1}w\psi : \mathbb{B}^k\to\mathbb{B}^k$ has a fixed point $\zeta\in\mathbb{B}^k$, by Brouwer's theorem. We get $w{\psi}(\zeta)={\psi}(\zeta)$ and by freeness of the action, this implies $w=1$.
\end{proof}

\begin{rem}
The first statement still holds if we replace ``CW structure'' by ``stratification''. More precisely, if $\mathcal{DV}$ admits a stratification such that for all $1\notin I\subset W$, the wall $Z_I$ is a (possibly empty) disjoint union of closed $|I|$-codimensional strata, if $\sigma\subset Z_I$ is a stratum and $w\in W\setminus\{1\}$ is such that $w\sigma\cap\mathcal{DV}\ne\emptyset$, then $w^{-1}\in I$ and $w\sigma$ is a stratum in $Z_w\cap Z_{wI\setminus\{1\}}$.

The following fact can be extracted from the end of the proof : if a discrete group acts freely on a regular CW complex, then the CW structure is equivariant if and only if the group acts on the cells. This is a consequence of Brouwer's fixed point theorem.
\end{rem}

We finish this section by illustrating the above results in the case of \emph{lens spaces}.

\begin{exemple}[Lens spaces]\label{lens_spaces}
Let $p,q\in\N^*$ be coprime integers, with $p>1$. Recall that we define a free action of $\Z/p\Z$ on $\Sph^3\subset \C^2$ by letting a generator of $\Z/p\Z$ act as
\[\rho(z_0,z_1):=(e^{2i\pi/p}z_0,e^{2iq\pi/p}z_1).\]
The quotient manifold $L(p;q):=\Sph^3/\left<\rho\right>$ is called a \emph{lens space}. We claim that the $\Z/p\Z$-equivariant CW structure on $\Sph^3$ inducing the usual CW structure on $L(p;q)$ \cite[Example 2.43]{hatcher} can be built from a Dirichlet--Voronoi domain.

To do this, we equip $\Sph^3$ with its natural round metric (which is obviously $\left<\rho\right>$-invariant) and recall that the geodesic distance to the ``central'' point $(1,0)$ is given by
\[\forall (z_0,z_1)\in\Sph^3,~d((1,0),(z_0,z_1))=\arccos\Re(z_0),\]
seeing again $\Sph^3$ inside $\C^2$. Then, the associated Dirichlet--Voronoi domain reads
\begin{align*}
\mathcal{DV}:=&\{x\in\Sph^3~;~d(1,x)\le d(1,\rho^kx),~\forall 0\le k<p\} \\
=&\{(z_0,z_1)\in\Sph^3~;~\arccos\Re(z_0)\le \arccos\Re(e^{2ik\pi/p}z_0),~\forall 0\le k<p\} \\
=&\{(z_0,z_1)\in\Sph^3~;~\Re((1-e^{2ik\pi/p})z_0)\ge0,~\forall 0\le k<p\}.
\end{align*}
The image $\pr_0(\mathcal{DV})$ of $\mathcal{DV}$ under the first projection $(z_0,z_1)\stackrel{\pr_0}\longmapsto z_0$ is represented in Figure \ref{angular}.

Observe that
\[\stackrel{\circ}{\mathcal{DV}}=\{(z_0,z_1)\in\Sph^3~;~\Re((1-e^{2ik\pi/p})z_0)>0,~\forall 0\le k<p\}\subset\{(z_0,z_1)~;~\Re(z_0)>0\}.\]
In particular, for $x\in~\stackrel{\circ}{\mathcal{DV}}$, we have $d(1,x)=\arccos(\Re(z_0))<\pi/2=\inj_1(\Sph^3)$. This means that $\stackrel{\circ}{\mathcal{DV}}~\subset B_{\Sph^3}(1,\inj\Sph^3)$ so that, by Proposition \ref{intDVcell}, the open domain $\stackrel{\circ}{\mathcal{DV}}$ is a 3-cell\footnote{This can be seen directly: if $(z_0,z_1)\in~\stackrel{\circ}{\mathcal{DV}}$, then $z_1$ lives in the open disk $B_{\C}(0,1)$ and, once $z_1$ is fixed, $z_0$ moves on an arc with radius $\sqrt{1-|z_1|^2}$ and angle in the interval $]-\tfrac{\pi}{p};\tfrac{\pi}{p}[$. Therefore, $\stackrel{\circ}{\mathcal{DV}}~\simeq B_{\C}(0,1)\times]-1;1[$ is a 3-cell.}.

\begin{center}
\begin{figure}
\begin{tikzpicture}[scale=2.2]
	\coordinate (z) at (0,0);
	\coordinate (p) at ({cos(180/7)},{sin(180/7)});
	\coordinate (pm) at ({cos(180/7)},{-sin(180/7)});
	\coordinate (a) at (1,0);
	\coordinate (t) at ({0.75*cos(360/36)},{0.75*sin(360/36)});
	
	\draw[very thick] (z)--(p) (z)--(pm);
	\draw[very thick] (pm) arc (-180/7:180/7:1);
	\draw[opacity=0.5] (z)--(a) (z)--(-1,0) (0,1)--(0,-1);
	\draw (z) circle (1cm);
	
	
	\draw (z) node[below left]{$0$};
	\draw (a) node[right]{$1$};
	\draw (p) node[above right]{$e^{i\pi/p}$};
	\draw (pm) node[below right]{$e^{-i\pi/p}$};
	\draw (t) node[above]{$z_0$};
	
	\fill[fill=black] (z) circle (.75pt);
	\fill[fill=black] (p) circle (.75pt);
	\fill[fill=black] (pm) circle (.75pt);
	\fill[fill=black] (a) circle (.75pt);
	\fill[fill=black] (t) circle (1pt);
	
	\begin{scope}[rotate=360/7]
	\coordinate (t) at ({0.75*cos(360/36)},{0.75*sin(360/36)});
	\fill[fill=black] (t) circle (1pt);
	\draw (t) node[below left]{$\rho z_0$};
	\end{scope}
	\begin{scope}[rotate=360*2/7]
	\coordinate (t) at ({0.75*cos(360/36)},{0.75*sin(360/36)});
	\fill[fill=black] (t) circle (1pt);
	\draw (t) node[below left]{$\rho^2 z_0$};
	\end{scope}
	\begin{scope}[rotate=360*3/7]
	\coordinate (t) at ({0.75*cos(360/36)},{0.75*sin(360/36)});
	\end{scope}
	\begin{scope}[rotate=360*4/7]
	\coordinate (t) at ({0.75*cos(360/36)},{0.75*sin(360/36)});
	\end{scope}
	\begin{scope}[rotate=360*5/7]
	\coordinate (t) at ({0.75*cos(360/36)},{0.75*sin(360/36)});
	\fill[fill=black] (t) circle (1pt);
	\draw (t) node[above left]{$\rho^{-2} z_0$};
	\end{scope}
	\begin{scope}[rotate=360*6/7]
	\coordinate (t) at ({0.75*cos(360/36)},{0.75*sin(360/36)});
	\fill[fill=black] (t) circle (1pt);
	\draw (t) node[above left]{$\rho^{-1} z_0$};
	\end{scope}
	
	\draw[thick, fill=black, opacity=0.3] (0,0)--(pm) arc(-180/7:180/7:1) -- cycle;
\end{tikzpicture}
\caption{The first projection $\pr_0(\mathcal{DV})$ of $\mathcal{DV}$ is the shaded angular sector.}\label{angular}
\end{figure}
\end{center}

We now look at the walls $Z_I$ of $\mathcal{DV}$, where $I\subset\Z/p\Z\setminus\{0\}$. When $I=\{k\}$ for some $0\ne k\in\Z/p\Z$, then it is immediate to check (e.g. using Figure \ref{angular}) that
\begin{equation}\label{2-walls}
\pr_0(Z_{k})=\pr_0(\mathcal{DV}\cap\rho^k\mathcal{DV})=\left\{\begin{array}{cc}[0;1]\cdot e^{\pm i\pi/p} & \text{if $k=\pm1$}, \\[.5em] \{0\} & \text{otherwise.} \end{array}\right.
\end{equation}
In particular, for $k\ne\pm1$, we have $Z_{k}=\{(0,z_1)\in\Sph^3\}\simeq\Sph^1$ is \emph{not} a disjoint union of 1-codimensional cells. Therefore, Proposition \ref{cellsfromDV} fails in this case\footnote{This happens because the distance on $\Sph^3$ ``doesn't see'' the $z_1$-factor, on which $\rho$ still acts non-trivially. Besides, observe that $Z_I\simeq\Sph^1$ as soon as $I\subset\Z/p\Z\setminus\{0\}$ has $|I|\ge2$.}. However, this circle turns out to be the unique obstruction; if we choose a $\Z/p\Z$-equivariant cell structure on it --- which is easy --- then we can obtain a CW structure on $\Sph^3$ using the walls $Z_1$ and $Z_{-1}=\rho^{-1}Z_1$.

Indeed, choosing for instance $z_1=1\in\Sph^1$, the orbit $\{\rho^k(0,z_1)~;~0\le k<p\}\subset\Sph^3$ forms the 0-cells, while the $p$ segments joining them give the 1-cells. Observe that we could have obtained a similar structure by taking again a Dirichlet--Voronoi domain for $\Z/p\Z$ acting on $\Sph^1$. Now, the wall
\[Z_1=\{(z_0,z_1)\in\Sph^3~;~e^{-i\pi/p}z_0\in[0;1]\}=\{(z_0,z_1)\in\Sph^3~;~z_0=e^{i\pi/p}\sqrt{1-|z_1|^2}\}\simeq\overline{B_\C(0,1)}\]
actually is a closed 2-cell and $Z_1\setminus\Sph^1=\{(z_0,z_1)\in Z_1~;~z_0\ne0\}=B_\C(0,1)$ is the representative 2-cell we look for. Completing with the 3-cell $\stackrel{\circ}{\mathcal{DV}}$, we obtain the stated CW decomposition of $\Sph^3$.

The situation can be visualized in the Figure \ref{vizu_lens}, and the resulting cellular chain complex for $L(p;q)$ is easily obtained:
\[\left(\xymatrix{\Z[\left<\rho\right>]\ar^{1-\rho^{-1}}[r] & \Z[\left<\rho\right>] \ar^{1+\rho+\cdots+\rho^{p-1}}[rr] & & \Z[\left<\rho\right>] \ar^{\rho-1}[r] & \Z[\left<\rho\right>]}\right)\otimes_{\Z[\left<\rho\right>]}\Z=\xymatrix{\Z\ar^0[r] & \Z\ar^p[r] & \Z\ar^0[r] & \Z.}\]

\begin{center}
\begin{figure}
\begin{tikzpicture}[scale=2]
	\coordinate (z) at (0,0,0);
	\coordinate (a) at (1,0,0);
	\coordinate (mam) at (0,1/2,0);
	\coordinate (ma) at (-1,0,0);
	\coordinate (mmam) at (0,-1/2,0);
	
	\draw plot [smooth,tension=1.5] coordinates {(a) (mam) (ma)};
	\draw plot [smooth,tension=1.5] coordinates {(ma) (mmam) (a)};
	
	\draw[opacity=0.5] plot [smooth,tension=1.5] coordinates {(a) (0,0.8,0) (ma)};
	\draw[opacity=0.5] plot [smooth,tension=1.5] coordinates {(ma) (0,-0.8,0) (a)};
	
	\draw[ultra thick] (-1,0,0) arc (180:360:1 and 0.2);
	\draw[ultra thin] (-1,0,0) arc (180:180+360/7:1 and 0.2) coordinate (R);
	\draw[ultra thin] (-1,0,0) arc (180:180+2*360/7:1 and 0.2) coordinate (S);
	\draw[ultra thin] (-1,0,0) arc (180:180+3*360/7:1 and 0.2) coordinate (T);
	\draw[ultra thick,dashed] (1,0,0) arc (0:180:1 and 0.2);
	
	\draw (1.2,0.5,0) node[right]{$\mathcal{DV}\cap\rho\mathcal{DV}$};
	\draw[->,>=stealth'] (1.2,0.5,0)--(0.65,0.4,0);
	\draw (1.2,-0.5,0) node[right]{$\mathcal{DV}\cap\rho^{-1}\mathcal{DV}$};
	\draw[->,>=stealth'] (1.2,-0.5,0)--(0.65,-0.4,0);
	
	\draw[opacity=0.5] (1.25,0.8,0) node[right]{$\rho\mathcal{DV}\cap\rho^{2}\mathcal{DV}$};
	\draw[->,>=stealth',opacity=0.5] (1.2,0.8,0)--(0.65,0.6,0);
	\draw[opacity=0.5] (1.25,-0.8,0) node[right]{$\rho^{-1}\mathcal{DV}\cap\rho^{-2}\mathcal{DV}$};
	\draw[->,>=stealth',opacity=0.5] (1.2,-0.8,0)--(0.65,-0.6,0);
	
	\draw (-1.3,0,0) node[left]{$\mathcal{DV}\cap\rho^\pm\mathcal{DV}\simeq\Sph^1$};
	\draw[->,>=stealth'] (-1.3,0,0)--(-1.05,0,0);
	
	\draw (z) node{$\mathcal{DV}$};
	\draw[opacity=0.5] (mam) node[above]{$\rho\mathcal{DV}$};
	\draw[opacity=0.5] (mmam) node[below]{$\rho^{-1}\mathcal{DV}$};
	
	\fill[fill=black] (R) circle (1pt);
	\fill[fill=black] (S) circle (1pt);
	\fill[fill=black] (T) circle (1pt);
	\draw[thick, dotted] plot [smooth,tension=1.45] coordinates{(R) (-1/2,1/6,-1/4) (mam)};
	\draw[thick, dotted] plot [smooth,tension=1.45] coordinates{(S) (1/18,1/6,-1/4) (mam)};
	\draw[thick, dotted] plot [smooth,tension=1.45] coordinates{(T) (1/2,1/6,-1/4) (mam)};
	\draw[thick, dotted] plot [smooth,tension=1.45] coordinates{(R) (-1/3.3,-1/3.7,1/4) (mmam)};
	\draw[thick, dotted] plot [smooth,tension=1.45] coordinates{(S) (1/4,-1/3.7,1/4+0.2*0) (mmam)};
	\draw[thick, dotted] plot [smooth,tension=1.45] coordinates{(T) (1/1.6,-1/3.7,1/4) (mmam)};
\end{tikzpicture}
\caption{The classical visualization of the lens space $L(p;q)$.}\label{vizu_lens}
\end{figure}
\end{center}

A recursion of the same procedure can be applied to higher dimensional lens spaces $L(p;q_1,\dotsc,q_n)$. Indeed, the open Dirichlet--Voronoi domain is still a $(2n+1)$-cell in $\Sph^{2n+1}$ (thanks to Proposition \ref{intDVcell}) and the two walls $Z_{\pm1}$ bounding it are again closed $2n$-cells, glued along their common boundary, a $\Z/p\Z$-stable sphere $\Sph^{2n-1}$. The resulting action of $\Z/p\Z$ on $\Sph^{2n-1}$ is exactly the same as the action associated to the lower lens space $L(p;q_2q_1^{-1},\dotsc,q_nq_1^{-1})$ ($q_1^{-1}$ denoting an inverse of $q_1$ mod $p$), twisted by the automorphism $\Z/p\Z\stackrel{\tiny{\times q_1}}\longto\Z/p\Z$. Therefore, we can consider a Dirichlet--Voronoi domain for this action, that leads to the stated $\Z/p\Z$-equivariant CW structure on $\Sph^{2n+1}$, and thus to the CW structure on $L(p;q_1,\dotsc,q_n)$, which is easily seen to be the same as in \cite{hatcher}.
\end{exemple}

\part{The case of flag manifolds: injectivity radius and a new equivariant CW structure on $O(3)/O(1)^3$}

In this second part, we aim to apply the above general results on Riemannian Dirichlet--Voronoi domains to the construction of equivariant cell structures on real flag manifolds, the acting group being the Weyl group $W$.

As we have seen, the injectivity radius plays an important role in this method. This is why the first section is devoted to the computation of this invariant for all complex and real flag manifolds, in terms of their root systems. As we will see, this allows us to recover the B\"{o}ttcher--Wenzel inequality \cite{bottcher-wenzel} for skew-hermitian matrices and the Bloch--Iserles inequalities \cite{bloch-iserles} for skew-symmetric (real) matrices.

Next, we focus on the particular case of the flag manifold $O(3)/O(1)^3$ of $SL_3(\R)$, for which we provide a new $W$-equivariant cell structure and compute the associated cellular chain complex. The proof is based on the quaternionic description of $O(3)/O(1)^3$ given in \cite[\S 4.4]{chirivi-garnier-spreafico} and is still to be adapted to the higher cases.

\section{The injectivity radius of flag manifolds}

\subsection{General setting and notation}\label{setting}
Throughout, we let $K$ be a connected semisimple compact Lie group, equipped with a maximal torus $T\le K$ and Weyl group $W=N_K(T)/T$. We also denote by $\Phi\subset\ho(T,\Sph^1)$ the root system of $(K,T)$, with chosen positive (resp. simple) roots set $\Phi^+$ (resp. $\Pi$). For a root $\alpha\in\ho(T,\Sph^1)$, we denote by ${\rm d}\alpha:={\rm d}_1\alpha : \mathfrak{t}\to i\R$ its differential; an element of $(i\mathfrak{t})^*$, where $\mathfrak{t}$ denotes the Lie algebra of $T$.

Let $G:=K^\C$ be the complexification of $K$, with maximal torus $T^\C$. This is a connected, semisimple complex algebraic group, whose Lie algebra $\mathfrak{g}$ is the complexification of the Lie algebra $\mathfrak{k}$ of $K$ and $\mathfrak{h}:=\mathfrak{t}\otimes\C$ is the Lie algebra of $T^\C$: a Cartan subalgebra of $\mathfrak{g}$. The \textit{flag manifold} $\mathcal{F}_K$ is the homogeneous space $K/T$ which, thanks to the Iwasawa decomposition, is diffeomorphic to $G/B$, where $B\le G$ is the Borel subgroup of $G$ containing $T$, associated to $\Phi^+$. This allows us to regard $\mathcal{F}_K$ as a smooth complex projective variety.

Since the group $T$ is abelian, the Weyl group $W$ acts freely on the right of $\mathcal{F}_K$ as follows: if $w\in W$ is represented by $\dot{w}\in N_K(T)$ and if $x\in\mathcal{F}_K$ is represented by $g\in K$, then $x\cdot w$ is defined as the class of $g\dot{w}$. It should be noted that the resulting action on the projective variety $G/B$ is \textit{not} algebraic.

On the other hand, the flag manifold $\mathcal{F}_K$ can be canonically endowed with the structure of a Riemannian homogeneous space as follows: let $\kappa : \mathfrak{g}\times\mathfrak{g}\to\C$ be the Killing form on $\mathfrak{g}$. It is non-degenerate since $\mathfrak{g}$ is (semi-)simple and its restriction to $\mathfrak{k}$ is negative-definite, as $K$ is compact \cite[Lemma 7.36]{besse_einstein}. Hence, the form $\left<\cdot,\cdot\right>:=-\kappa$ is an inner product on $\mathfrak{k}=T_1K$, inducing a bi-invariant Riemannian metric $\widetilde{g}_\kappa$ on $K$. This metric being bi-invariant, it descends to a $K$-invariant Riemannian metric $g_\kappa$ on $\mathcal{F}_K=K/T$, called the \textit{standard normal homogeneous} metric. In particular, the free action of the Weyl group $W$ on $\mathcal{F}_K$ is isometric. Moreover, by \cite[Theorem 2.5]{ye-wong-lim}, the geodesics of $(\mathcal{F}_K,g_K)$ are orbits of one-parameter subgroups of $K$.

Observe moreover that if $K$ is simple and if $K_{\rm sc}$ (resp. $K_{\rm ad}$) denotes the simply-connected (resp. adjoint) group of the same type as $K$, then there are isogenies $K_{\rm sc}\twoheadrightarrow K\twoheadrightarrow K_{\rm ad}$ whose differentials are the identity on $\mathfrak{k}$. Moreover, these isogenies induce isomorphisms of flag manifolds $K_{\rm sc}/T_{\rm sc}\stackrel{\tiny{\sim}}\to K/T\stackrel{\tiny{\sim}}\to K_{\rm ad}/T_{\rm ad}$ that are locally isometric. Therefore, these are $W$-equivariant isometries and we may talk about \textit{the} (standard homogeneous) flag manifold $\mathcal{F}_\Phi:=\mathcal{F}_K$, associated to the root system $\Phi$ of $K$.

Recall the root space decomposition 
\[\mathfrak{g}=\mathfrak{h}\oplus\bigoplus_{\alpha\in\Phi}\mathfrak{g}_\alpha,\]
where $\mathfrak{g}_\alpha=\{x\in\mathfrak{g}~|~[h,x]={\rm d}\alpha(h)x,~\forall h\in\mathfrak{h}\}$. Consider the Chevalley--Serre generators $(e_\alpha,f_\alpha,h_\alpha)_{\alpha\in\Phi^+}$ of $\mathfrak{g}$. Then for $\alpha\in\Phi^+$, we have $\mathfrak{g}_{\alpha}=\C e_\alpha$ and $\mathfrak{g}_{-\alpha}=\C f_\alpha$ and the root decomposition reads
\[\mathfrak{g}=\bigoplus_{\alpha\in\Pi}\C h_\alpha\oplus\bigoplus_{\alpha\in\Phi^+}(\C e_\alpha\oplus\C f_\alpha).\]

Now, for $\alpha\in\Phi^+$, we let
\[u_\alpha:=\frac{e_\alpha-f_\alpha}{2},~v_\alpha:=\frac{i(e_\alpha+f_\alpha)}{2},~w_\alpha:=\frac{ih_\alpha}{2}.\]
Then, the collection $(u_\alpha,v_\alpha,w_\alpha)_{\alpha\in\Phi^+}$ generates $\mathfrak{k}$, in the sense that
\[\mathfrak{k}=\mathfrak{t}\oplus\bigoplus_{\alpha\in\Phi^+}(\R u_\alpha\oplus\R v_\alpha)=\bigoplus_{\alpha\in\Pi}\R w_\alpha\oplus\bigoplus_{\alpha\in\Phi^+}(\R u_\alpha\oplus\R v_\alpha).\]
We have the \textit{quaternionic} relations
\[\forall \alpha\in\Phi^+,~[u_\alpha,v_\alpha]=w_\alpha,~[v_\alpha,w_\alpha]=u_\alpha,~[w_\alpha,u_\alpha]=v_\alpha,\]
that is, the triple $(u_\alpha,v_\alpha,w_\alpha)$ generates a subalgebra of $\mathfrak{k}$ isomorphic to $\mathfrak{so}_3(\R)$. It may also be observed that the family $\{u_\alpha,v_\alpha\}_{\alpha\in\Phi^+}$ is orthogonal (for the product $\left<\cdot,\cdot\right>=-\kappa$) and that this family is itself orthogonal to each $w_\alpha$ (although the family $\{w_\alpha\}_{\alpha\in\Pi}$ is not orthogonal in general). In particular, the orthogonal of the toral subalgebra $\mathfrak{t}$ is the subspace $\mathfrak{p}:=\bigoplus^{\perp}_{\alpha\in\Phi^+}(\R u_\alpha \oplus^\perp\R v_\alpha)$. Moreover, if $t_\alpha\in\mathfrak{h}$ denotes the dual of $\alpha\in\Phi$, (i.e. $\kappa(t_\alpha,-)={\rm d}\alpha$, the existence being guaranteed by non-degeneracy of $\kappa$), then we have 
\[\left<u_\alpha,u_\alpha\right>=-\kappa(u_\alpha,u_\alpha)=\frac{-\kappa(e_\alpha-f_\alpha,e_\alpha-f_\alpha)}{4}=\frac{\kappa(e_\alpha,f_\alpha)}{2}=\frac{1}{\kappa(t_\alpha,t_\alpha)}=\frac{1}{\|\alpha\|^2},\]
with $\|\alpha\|^2:=\kappa(t_\alpha,t_\alpha)={\rm d}\alpha(t_\alpha)>0$. In the same fashion, we see that $\|v_\alpha\|=\|w_\alpha\|=1/\|\alpha\|=\|u_\alpha\|$. These considerations yield an explicit formula for the norm on $\mathfrak{k}$: writing an element $X\in\mathfrak{k}$ as $X=X^\mathfrak{t}+X^\mathfrak{p}$, where $X^\mathfrak{p}=\sum_{\alpha\in\Phi^+}\lambda_\alpha u_\alpha+\mu_\alpha v_\alpha$ (for $\lambda_\alpha,\mu_\alpha\in\R$), we have
\[\|X\|^2=\|X^\mathfrak{t}\|^2+\|X^\mathfrak{p}\|^2=\|X^\mathfrak{t}\|^2+\sum_{\alpha\in\Phi^+}\frac{\lambda_\alpha^2+\mu_\alpha^2}{\|\alpha\|^2}=\sum_{\alpha\in\Phi}|{\rm d}\alpha(X^\mathfrak{t})|^2+\sum_{\alpha\in\Phi^+}\frac{\lambda_\alpha^2+\mu_\alpha^2}{\|\alpha\|^2}\]
so that we obtain
\begin{equation}\label{norm}\tag{N}
\|X\|^2=\sum_{\alpha\in\Phi^+}\left(2|{\rm d}\alpha(X^\mathfrak{t})|^2+\frac{\lambda_\alpha^2+\mu_\alpha^2}{\|\alpha\|^2}\right).
\end{equation}
This will be useful when we estimate the norm of the map $\mathfrak{k}\stackrel{\tiny{{\rm ad}}}\longto\mathfrak{gl}(\mathfrak{k})$.

\subsection{Lemmata and main result}
From now on, we assume that the root system $\Phi$ is irreducible.

Recall the Coxeter number $h$ and its dual $h^\vee$ are defined as $h=1+\sum_{\alpha\in\Pi}n_\alpha$ and $h^\vee=1+\sum_{\alpha\in\Pi}n_\alpha^\vee$, where $\alpha_0=\sum_{\alpha\in\Pi}n_\alpha\alpha$ and $\alpha_0^\vee=\sum_{\alpha\in\Pi}n_\alpha^\vee\alpha^\vee$ are the highest long root and highest short coroot of $\Phi$, respectively. If $\rho:=\tfrac12\sum_{\alpha\in\Phi^+}\alpha$ denotes the Weyl vector and $\rho^\vee:=\tfrac12\sum_{\alpha\in\Phi^+}\alpha^\vee$, then we have $h=1+\left<\rho^\vee,\alpha_0\right>$ and $h^\vee=1+\left<\alpha_0^\vee,\rho\right>$. While standard, the following fact will be useful for our analysis, so we include an independent proof for completeness.

\begin{lem}[{\cite[Lemma 4]{suter}}]\label{hvee}
If $\mathfrak{g}$ is simple, the dual Coxeter number $h^\vee=1+\left<\alpha_0^\vee,\rho\right>$ satisfies
\[h^\vee=\frac{1}{\|\alpha_0\|^2}.\] 
\end{lem}
\begin{proof} For each $\alpha\in\Phi^+$, recall the root vectors $e_\alpha$ and $f_\alpha$ and let $e_{-\alpha}:=\frac{\|\alpha\|^2}{2}f_\alpha$ (so that $\kappa(e_\alpha,e_{-\alpha})=1$). We introduce the \textit{Casimir element} for the adjoint representation: \[c:=c_{\rm ad}=\sum_{\alpha\in\Pi}{\rm ad}(h_\alpha){\rm ad}(h_\alpha^\perp)+\sum_{\alpha\in\Phi^+}{\rm ad}(e_\alpha){\rm ad}(e_{-\alpha})+{\rm ad}(e_{-\alpha}){\rm ad}(e_\alpha)\in\mathfrak{gl}(\mathfrak{g}).\] Observe that if $\alpha,\beta\in\Pi$, then we have
\[\sum_{\delta\in\Pi}\alpha(h_\delta)\beta(h_\delta^\perp)=\sum_{\delta\in\Pi}\frac{\|\beta\|^2}{2}\left<\delta^\vee,\alpha\right>\left<h_\beta,h_\delta^\perp\right>=\frac{\|\beta\|^2}{2}\left<\beta^\vee,\alpha\right>=\alpha(t_\beta)=\kappa(\alpha,\beta)\]
and thus, \[\forall U,V\in\mathfrak{t},~\left<U,V\right>=\sum_{\alpha\in\Pi}U(h_\alpha)V(h_\alpha^\perp).\] Now, since $e_{\alpha_0}$ is a highest weight vector for the adjoint representation, ${\rm ad}(e_\alpha)(e_{\alpha_0})=[e_\alpha,e_{\alpha_0}]=0$ for all $\alpha\in\Phi^+$. Since $c={\rm id}_{\mathfrak{g}}$, we find \begin{align*}
e_{\alpha_0}&=c(e_{\alpha_0})=\sum_{\alpha\in\Pi}{\rm ad}(h_\alpha){\rm ad}(h_\alpha^\perp)(e_{\alpha_0})+\sum_{\alpha\in\Phi^+}{\rm ad}(e_\alpha){\rm ad}(e_{-\alpha})(e_{\alpha_0}) \\
&=\sum_{\alpha\in\Pi}\alpha_0(h_\alpha^\perp){\rm ad}(h_\alpha)(e_{\alpha_0})+\sum_{\alpha\in\Phi^+}{\rm ad}[e_\alpha,e_{-\alpha}](e_{\alpha_0}) \\
&=\sum_{\alpha\in\Pi}\alpha_0(h_\alpha)\alpha_0(h_\alpha^\perp)e_{\alpha_0}+\sum_{\alpha\in\Phi^+}\frac{\|\alpha\|^2}{2}[h_\alpha,e_{\alpha_0}]=\left<\alpha_0,\alpha_0\right>e_{\alpha_0}+\sum_{\alpha\in\Phi^+}\frac{\|\alpha\|^2}{2}\alpha_0(h_\alpha)e_{\alpha_0} \\
&=\left(\|\alpha_0\|^2+\sum_{\alpha\in\Phi^+}\frac{\|\alpha\|^2}{2}\left<\alpha^\vee,\alpha_0\right>\right)e_{\alpha_0}=\left(\|\alpha_0\|^2+\sum_{\alpha\in\Phi^+}\left<\alpha,\alpha_0\right>\right)e_{\alpha_0} \\
&=\left<\alpha_0,\alpha_0+2\rho\right>e_{\alpha_0}
\end{align*} and thus, \[1=\|\alpha_0\|^2\left(1+\left<\frac{2\alpha_0}{\|\alpha_0\|^2},\rho\right>\right)=\|\alpha_0\|^2(1+\left<\alpha_0^\vee,\rho\right>)=\|\alpha_0\|^2h^\vee.\]
\end{proof}

We start with the following key lemma on the minimal length a closed geodesic in $\mathcal{F}_\Phi$ can have.
\begin{lem}\label{mingeo}
The minimal length of a non-trivial closed geodesic in $\mathcal{F}_\Phi$ is $2\pi\sqrt{h^\vee}$.
\end{lem}

\begin{proof}
As mentioned in \S \ref{setting}, we may assume that $\mathcal{F}_\Phi=K/T$, where $K=K_{\rm sc}$ is the simply-connected (simple) compact Lie group with root system $\Phi$. A closed geodesic $\gamma : [0,1]\to\mathcal{F}_K$ sends $s\in[0,1]$ to $\gamma(s)=e^{sX}\cdot T$, for some $X\in\mathfrak{k}$ with $\|X\|=\leng(\gamma)$. Saying that $\gamma$ is non-trivial amounts to say that $X\notin\mathfrak{t}$ and  since $\gamma$ is closed, we have $e^X\in T$. We first prove that $\|X\|\ge2\pi\sqrt{h^\vee}$.

We choose an element $g\in K$ such that ${\rm Ad}_g(X)=:Z\in\mathfrak{t}$ and assume for contradiction that $2\pi\sqrt{h^\vee}=2\pi/\|\alpha_0\|>\|X\|=\|Z\|$. We claim the equality of centralizers
\[C_G(e^Z)=C_G(Z).\]
The inclusion $C_G(Z)\le C_G(e^Z)$ is trivial. Since $K$ is simply-connected, $G=K^\C$ is simply-connected as well, so that by \cite[Theorem 2.11]{hum_conj}, the centralizer $C_G(e^Z)=C_G(e^Z)^o$ is connected and using \cite[Theorem 2.2]{hum_conj}, it is generated by $T^\C$ and those root subgroups $U_\alpha$ ($\alpha\in\Phi^+$) for which $\alpha(e^Z)=1$. Thus, it suffices to prove that $U_\alpha\le C_G(Z)$ if $\alpha(e^Z)=1$. Take $\alpha\in\Phi^+$ such that $1=\alpha(e^Z)=e^{{\rm d}\alpha(Z)}$. There is some $k\in\Z$ such that ${\rm d}\alpha(Z)=2ik\pi$ and by the Cauchy--Schwarz inequality, we arrive at
\[2|k|\pi=|{\rm d}\alpha(Z)|=\left|\left<t_\alpha,Z\right>\right|\le\|t_\alpha\|\|Z\|=\|\alpha\|\|Z\|\le\|\alpha_0\|\|Z\|<2\pi,\]
so that $k=0$ and ${\rm d}\alpha(Z)=0$. This ensures that $[Z,X_\alpha]=0$ for each $X_\alpha\in\mathfrak{g}_\alpha$ so that $U_\alpha$ centralizes $Z$, as claimed. Now, because the elements $e^X$ and $e^Z$ lie in $T$ and are $G$-conjugate, we may apply \cite[Lemma 4.33]{adams_lie} so find some $\dot{w}\in N_G(T^\C)$ such that $e^X={\rm Ad}_{g^{-1}}(e^Z)={\rm Ad}_{\dot{w}}(e^Z)$. This yields $g\dot{w}\in C_G(e^Z)=C_G(Z)$ and thus, $X={\rm Ad}_{g^{-1}}(Z)={\rm Ad}_{\dot{w}}(Z)\in\mathfrak{t}$, an absurdity.

Conversely, we find some $X\in\mathfrak{k}\setminus\mathfrak{t}$ with norm $\|X\|=2\pi\sqrt{h^\vee}$ and such that $e^X\in T$. To the $\mathfrak{sl}_2$-triple $(e_{\alpha_0},f_{\alpha_0},h_{\alpha_0})$ corresponds an inclusion $\mathfrak{sl}_2(\C)\hookrightarrow\mathfrak{g}$, which is the differential of a map $\phi_0 : SL_2(\C)\hookrightarrow G$ (restricting to an inclusion $\widetilde{\phi}_0 : SU(2)\hookrightarrow K$). Then, the reflection $s_{\alpha_0}\in W$ may be represented by the element $\dot{s}_{\alpha_0}=\phi_0\left(\begin{smallmatrix}0 & 1 \\ -1 & 0\end{smallmatrix}\right)$ and we have
\begin{align*}
\dot{s}_{\alpha_0}&=\phi_0\begin{pmatrix}0 & 1 \\ -1 & 0\end{pmatrix}=\phi_0\left(\exp\left(\frac{\pi}{2}\begin{pmatrix}0 & 1 \\ -1 & 0\end{pmatrix}\right)\right) \\
&=\exp\left(\frac{\pi}{2}{\rm d}\phi_0\begin{pmatrix}0 & 1 \\ -1 & 0\end{pmatrix}\right)=\exp\left(\frac{\pi(e_{\alpha_0}-f_{\alpha_0})}{2}\right) \\
&=\exp(\pi u_{\alpha_0})
\end{align*}
so that if we let $X_0:=2\pi u_{\alpha_0}\notin\mathfrak{t}$, then $e^{X_0}=\dot{s}^2_{\alpha_0}\in T$ and $\|X_0\|=2\pi\|u_{\alpha_0}\|=\frac{2\pi}{\|\alpha_0\|}$.
\end{proof}

We now give a sharp estimate on the sectional curvature of flag manifolds:
\begin{lem}\label{curvat}
Let $\mathcal{K}$ be the sectional curvature of $\mathcal{F}_\Phi$. For any plane $P\le T_1\mathcal{F}_\Phi$, we have
\[0\le\mathcal{K}(P)\le\frac{1}{h^\vee}.\]
Moreover, this estimate is sharp: if $P_0:=\R u_{\alpha_0}\oplus\R v_{\alpha_0}$, then $\mathcal{K}(P_0)=1/h^\vee$.
\end{lem}

\begin{proof}
Let $(X,Y)$ be an orthonormal basis of $P\le T_1\mathcal{F}_\Phi=\mathfrak{p}$. By the O'Neill formula \cite[Theorems 3.61 \& 3.65]{gallot-hulin-lafontaine}, we have
\[0\le\mathcal{K}(P)=\mathcal{K}(X,Y)=\frac14\|[X,Y]^\mathfrak{p}\|^2+\|[X,Y]^\mathfrak{t}\|^2=\frac14\|[X,Y]\|^2+\frac34\|[X,Y]^\mathfrak{t}\|^2\le\|[X,Y]\|^2.\]
Since ${\rm Ad}_k\in\mathfrak{gl}(\mathfrak{k})$ is an isometry for all $k\in K$, to estimate the squared norm $\|[X,Y]\|^2$, we may assume that $Y=:H\in\mathfrak{t}$, in which case $X$ can be any element in the unit sphere of $\mathfrak{k}$. But writing $X=X^{\mathfrak{t}}+X^{\mathfrak{p}}$, we have $[X,H]=[X^{\mathfrak{p}},H]$, so that we may further assume that $X=X^{\mathfrak{p}}\in\mathfrak{p}$, in which case $[X,H]\in\mathfrak{p}$.

Therefore, it suffices to estimate the squared norm $\|[H,X]\|^2$ for $H\in\mathfrak{t}$ in the unit sphere and $X\in\mathfrak{p}$ in the (closed) unit ball. Following \S \ref{setting}, write $X=\sum_{\alpha\in\Phi^+}\lambda_\alpha u_\alpha+\mu_\alpha v_\alpha$ (where $\lambda_\alpha,\mu_\alpha\in\R$), so that by \eqref{norm}, we have
\[\|X\|^2=\sum_{\alpha\in\Phi^+}\frac{\lambda_\alpha^2+\mu_\alpha^2}{\|\alpha\|^2}\le1.\]
It is immediate to check that, for every $\alpha\in\Phi^+$, we have $[H,u_\alpha]=-i{\rm d}\alpha(H)v_\alpha$ and $[H,u_\alpha]=i{\rm d}\alpha(Z)u_\alpha$, so that
\[[H,X]=\sum_{\alpha\in\Phi^+}\lambda_\alpha[H,u_\alpha]+\mu_\alpha[H,v_\alpha]=i\sum_{\alpha\in\Phi^+}{\rm d}\alpha(H)(\mu_\alpha u_\alpha-\lambda_\alpha v_\alpha)\]
and, using \eqref{norm} again,
\begin{align*}
\|[H,X]\|^2&=\sum_{\alpha\in\Phi^+}|{\rm d}\alpha(H)|^2\frac{\lambda_\alpha^2+\mu_\alpha^2}{\|\alpha\|^2}\le\max_{\alpha\in\Phi}|{\rm d}\alpha(H)|^2\sum_{\beta\in\Phi^+}\frac{\lambda_\beta^2+\mu_\beta^2}{\|\beta\|^2} \\
&\le\max_{\alpha\in\Phi^+}|{\rm d}\alpha(H)|^2\le\|H\|^2\max_{\alpha\in\Phi^+}\|\alpha\|^2=\|\alpha_0\|^2 \\
&={1}/{h^\vee}.
\end{align*}

Conversely, the family $(\|\alpha_0\|u_{\alpha_0},\|\alpha_0\|v_{\alpha_0})$ is an orthonormal basis of $P_0$ and we have
\[\mathcal{K}(P_0)=\frac{\|\alpha_0\|^4}{4}\|\underbrace{[u_{\alpha_0},v_{\alpha_0}]^\mathfrak{p}}_{=0}\|^2+\|\alpha_0\|^4\|[u_{\alpha_0},v_{\alpha_0}]^\mathfrak{t}\|^2=\|\alpha_0\|^4\|w_{\alpha_0}\|^2=\|\alpha_0\|^2=1/h^\vee.\]
\end{proof}

Now, the Klingenberg Lemma \cite[Chapter 13, Proposition 2.13]{do_carmo_riemannian} (see also \cite[Proposition 2.6.8]{klingenberg_riemann}) tells us that, if $\ell=2\pi\sqrt{h^\vee}$ denotes the minimal length of a closed geodesic in $\mathcal{F}_\Phi$ (see Lemma \ref{mingeo}), then
\[\inj(\mathcal{F}_\Phi)\ge\min\left\{\frac{\pi}{\sqrt{\sup_P\mathcal{K}(P)}},\frac{\ell}{2}\right\}=\pi\sqrt{h^\vee}.\]
Since the reverse inequality $\inj(\mathcal{F}_\Phi)\le\ell/2=\pi\sqrt{h^\vee}$ is obvious, we arrive at the main result of this section:
\begin{theo}\label{injrad}
The injectivity radius of a standard homogeneous flag manifold $\mathcal{F}_\Phi$ is given by
\[\inj(\mathcal{F}_\Phi)=\pi\sqrt{h^\vee},\]
where $h^\vee$ is the dual Coxeter number of the (simple) root system $\Phi$. Moreover, this is the distance between 1 and any reflection of $W$, whose associated root is long.
\end{theo}

\begin{rem}\label{BWAC}
We make the following observations:
\begin{itemize}
\item In the proof of Lemma \ref{curvat}, we have seen that $\|[X,Y]\|^2\le\tfrac{1}{h^\vee}\|X\|^2\|Y\|^2$ for any two vectors $X,Y\in\mathfrak{k}$, where $\|\cdot\|$ is the Killing norm. In the case where $\mathfrak{g}$ is classical, we may interpret this in terms of the Frobenius norm $\|\cdot\|_F=\tfrac{1}{\lambda}\|\cdot\|$ by writing that 
\[\|[X,Y]\|_F\le\frac{\lambda}{\sqrt{h^\vee}}\|X\|_F\|Y\|_F.\]
The value of $\lambda$, depending on the type of the classical algebra, can be found in \cite[\S 12.6]{erdmann-wildon}. For instance, if $\mathfrak{k}=\mathfrak{su}(n)$, then $\lambda=\sqrt{2n}$ and we obtain that
\[\forall X,Y\in\mathfrak{su}(n),~\|[X,Y]\|_F\le\sqrt{2}\|X\|_F\|Y\|_F.\]
This is precisely the B\"{o}ttcher--Wenzel inequality \cite[Theorem 2.2]{bottcher-wenzel} for $\mathfrak{su}(n)$. This inequality turns out to hold on the bigger algebra $\mathfrak{gl}_n(\C)$ and is sharp. In type $C_n$, we have $\lambda=\sqrt{2(n+1)}$, so that we retrieve the same B\"{o}ttcher--Wenzel inequality for $\mathfrak{k}=\mathfrak{sp}(n)=\mathfrak{sp}_{n}(\C)\cap \mathfrak{u}(2n)$.

In types $B$ and $D$, we have $\mathfrak{k}=\mathfrak{so}_n(\R)$ and $\lambda^2=2\left\lfloor\frac{n-1}{2}\right\rfloor+\frac{(-1)^n-1}{2}=h^\vee$, so that we obtain the following improvement of the BW inequality:
\[\forall n\ge 5,~\sup_{X,Y\in\mathfrak{so}_n(\R)\setminus\{0\}}\frac{\|[X,Y]\|_F}{\|X\|_F\|Y\|_F}=1.\]
Moreover, as $\mathfrak{so}_3(\R)\simeq\mathfrak{su}(2)$, we have
\[\sup_{0\ne X,Y\in\mathfrak{so}_3(\R)}\frac{\|[X,Y]\|_F}{\|X\|_F\|Y\|_F}=\tfrac{1}{\sqrt{2}}.\]
Indeed, the map 
\[\left(\begin{smallmatrix}0 & x & z \\ -x & 0 & y \\ -z & -y & 0\end{smallmatrix}\right)\mapsto\tfrac12\left(\begin{smallmatrix}-iz & -ix-y \\ -ix+y & iz\end{smallmatrix}\right)\]
defines an isometry $(\mathfrak{so}_3(\R),\kappa)\stackrel{\tiny{\sim}}\to(\mathfrak{su}(2),\kappa)$ of Lie algebras and thus
\[\hspace{13.5mm}\sup_{0\ne X,Y\in\mathfrak{so}_3(\R)}\frac{\|[X,Y]\|_F}{\|X\|_F\|Y\|_F}=\sup_{0\ne X,Y\in\mathfrak{so}_3(\R)}\frac{\|[X,Y]\|_\kappa}{\|X\|_\kappa\|Y\|_\kappa}=\sup_{0\ne X,Y\in\mathfrak{su}(2)}\frac{\|[X,Y]\|_\kappa}{\|X\|_\kappa\|Y\|_\kappa}=\frac{1}{\sqrt{2}}.\]
These results were first derived in \cite[Theorem 6]{bloch-iserles}, with totally different methods. However, it should be noted that the authors also establish the displayed formula for $n=4$, which is unreachable by our method, as the Lie algebra $\mathfrak{so}_4(\R)\simeq\mathfrak{su}(2)\oplus\mathfrak{su}(2)$ is not simple.
\item It is proved in \cite[Theorem 2.11]{agrachev-chambrion} that the radius of $\mathcal{F}_{A_{n}}$ is given by
\[{\rm diam}(\mathcal{F}_{A_{n}})=\frac{\pi}{3}\sqrt{6n(n+2)}\]
and is the distance from $1$ to any Coxeter element of the Weyl group $W=\Sym_{n+1}$.
\item Our result for the injectivity radius in types $A_2$ and $C_3$ agrees with \cite[\S 3]{puttmann}.
\item The injectivity radius of the irreducible symmetric spaces was determined in \cite{yang1} and \cite{yang2}. Before that, the curvature and minimal length of closed geodesics in simply-connected symmetric spaces were given in \cite[Theorems 11.1 and 11.2]{helgason} and the results are very similar to the Lemmata \ref{mingeo} and \ref{curvat} above.

However, one should notice that, except for $A_1$, the flag manifolds are \emph{not} symmetric spaces. Otherwise, the decomposition $\mathfrak{k}=\mathfrak{t}\oplus\mathfrak{p}$ of the Lie algebra of $K$ should satisfy $[\mathfrak{t},\mathfrak{t}]\subset\mathfrak{t}$, $[\mathfrak{t},\mathfrak{p}]\subset\mathfrak{p}$ and $[\mathfrak{p},\mathfrak{p}]\subset\mathfrak{t}$ \cite[Lemma 7.69]{besse_einstein}, but the latter inclusion fails.

Finally, from \cite[Table 4.2]{yang1}, we see that the injectivity radius of the flag manifold associated to a simply-connected simple compact Lie group is half the injectivity radius of the group.
\end{itemize}
\end{rem}

\section{The case of $\mathcal{F}_{A_2}(\R)=O(3)/O(1)^3$}\label{O(3)/O(1)3}

Let $n\ge2$ and denote by $\mathcal{F}_n$ the flag manifold of $SL_n(\C)$ :
\[\mathcal{F}_n:=\mathcal{F}_{A_{n-1}}=SL_n(\C)/B\simeq SU(n)/T,\]
where $B\subset SL_n(\C)$ is the Borel subgroup of upper-triangular matrices and $T\subset SU(n)$ is the subgroup of diagonal matrices. The diffeomorphism between the two homogeneous spaces above is given by the Gram--Schmidt process: a particular case of Iwasawa decomposition. The split real form $SL_n(\R)$ of $SL_n(\C)$ endows the projective variety $\mathcal{F}_n$ with a real structure whose real points are given by
\[\mathcal{F}_n(\R)=SL_n(\R)/(B\cap SL_n(\R))\simeq SO(n)/S(O(1)^n)\simeq O(n)/O(1)^n,\]
where $S(O(1)^n)\subset SO(n)$ (resp. $O(1)^n\subset O(n)$) is the subgroup of diagonal matrices, a finite group isomorphic to $(\Z/2\Z)^{n-1}$ (resp. to $(\Z/2\Z)^n$). Here again, the diffeomorphism is induced by the Gram--Schmidt process. We denote by $g_n=g_\kappa$ the natural homogeneous metric on $\mathcal{F}_n$, induced by the Killing form $\kappa=2n\tr$ on $\mathfrak{su}(n)$.

The Weyl group $W:=N_{SU(n)}(T)/T\simeq\Sym_n$ is the symmetric group on $n$ letters and recall its action on $\mathcal{F}_n$, described at the beginning of \S \ref{setting}. Viewing $\mathcal{F}_n$ as the set of decompositions of $\C^n$ as a direct sum of pairwise orthogonal lines, this action coincides with the natural permutation action of $W=\Sym_n$ on the lines. 

\begin{rem}\label{fubini-study}
There is another natural metric on $\mathcal{F}_n$ induced by a power of the Fubini--Study metric. Indeed, we have just seen that there is an embedding 
\[\iota : \mathcal{F}_n\longhookrightarrow (\C\Pro^{n-1})^n\]
and endowing the space $(\C\Pro^{n-1})^n$ with the natural product metric $g_{\rm FS}$ of the Fubini--Study metric on each copy of $\C\Pro^{n-1}$, it is easy to see that we have
\[g_n=2n\iota^*(g_{\rm FS})\]
and in particular, $d_{\rm FS}\le d_n$. Moreover, $\inj((\C\Pro^{n-1})^n,g_{\rm FS})=\pi/2$ and for a matrix $A=(a_{i,j})$ in $SU(n)$, we have $d_{\rm FS}(\overline{1},\overline{A})^2=\sum_i\arccos(|a_{i,i}|)^2$. 
\end{rem}

We can give an elementary description of the Dirichlet--Voronoi domain in the simplest type where $n=2$, justifying in particular that the method of the first part can only be used for \textit{real} flag manifolds. This is done in the following example:
\begin{exemple}\label{examSU2}
For $K=SU(2)$, the normal homogeneous flag manifold $\mathcal{F}_2=SU(2)/T\simeq\Sph^2$ is the round 2-sphere. The non-trivial element $s$ of the Weyl group $W\simeq C_2$ acts as the antipode on $\Sph^2$, so that $\mathcal{DV}$ is a closed half-sphere. Therefore, its boundary $\partial\mathcal{DV}=Z_s=H_s\simeq\Sph^1$ is a circle, which is not a disjoint union of 1-cells.

However, the real locus $\mathcal{F}_2(\R)\simeq\Sph^1$ may be interpreted as a circle orthogonal to $H_s$, so that $\mathcal{DV}\cap\mathcal{F}_2(\R)$ is a closed half-circle, whose boundary $H_s\cap\mathcal{F}_2(\R)$ is a union of two (antipodal) points. Thus, we indeed obtain a $W$-equivariant cell structure on $\mathcal{F}_2(\R)$.

To work the case of $\mathcal{F}_2(\C)$ out, we have to take a Dirichlet--Voronoi domain for ${\rm Stab}_W(H_s)=W$ acting on $H_s$, centered at some \emph{chosen} point in $H_s$. This gives an equivariant cell structure on $H_s$ and, adding the top cells $\interior{\mathcal{DV}}$ and $\interior{\mathcal{DV}}\cdot s$, we obtain an equivariant cell decomposition of $\mathcal{F}_2$. This is the trivial decomposition mentioned in \cite[\S 4.4]{chirivi-garnier-spreafico}.
\end{exemple}

\subsection{Bi-invariance of the quaternionic metric on $\mathcal{F}_{A_2}(\R)$}\label{biinvariance}

The 3-sphere $\Sph^3$ is a compact Lie group, when viewed as the unit quaternions \cite[Chap. 5]{URG}, and is isomorphic to $SU(2)$. Among its well-known finite subgroups \cite[Theorem 5.12]{URG}, it has the \textit{binary octahedral group} $\mathcal{O}$ of order 48, itself admitting the \textit{quaternion group} $\mathcal{Q}_8$ of order 8 as a normal subgroup, with factor group $\mathcal{O}/\mathcal{Q}_8\simeq\Sym_3$. As observed in \cite[\S 4.4]{chirivi-garnier-spreafico}, the universal cover $\Sph^3\stackrel{\tiny{\mu}}\longtwoheadrightarrow SO(3)$ induces a tower of coverings:
\[\xymatrix@R-1pc{\Sph^3 \ar@/_10pt/@{->>}_{\mathcal{Q}_8}[rrd] \ar@{->>}^<<<<<{\{\pm1\}}[r] & SO(3) \ar@{->>}^<<<<<{\{\pm1\}^2}[r] & SO(3)/\{\pm1\}^2 \eq[d] \ar@{->>}^<<<<<{\Sym_3}[r] & \Sph^3/\mathcal{O} \\ & & \mathcal{F}_3(\R) \ar@/_10pt/@{->>}_{\Sym_3}[ru] & }\]
where a group near an arrow denotes the fiber of the map. The resulting action of $\Sym_3$ on $\mathcal{F}_3(\R)=\Sph^3/\mathcal{Q}_8$ coincides with the action of the Weyl group. In other words, the projection map
\begin{equation}\label{projection_map}
\Pi : \Sph^3\longtwoheadrightarrow\mathcal{F}_3(\R)
\end{equation}
intertwines the actions of $\mathcal{O}$ and of the Weyl group $W=\Sym_3$. 

Moreover, the standard round metric on $\Sph^3$, being bi-invariant, descends to a bi-invariant metric on $SO(3)$, itself inducing a normal homogeneous metric $g_8$ on $\mathcal{F}_3(\R)$. This metric must therefore be a multiple of the metric $g_\kappa$ induced by the restriction to $\mathfrak{so}_3(\R)$ of the Killing form $\kappa=6\tr$ on the compact Lie algebra $\mathfrak{su}(3)$. The proportionality constant can be determined using the diameters of the metrics. Indeed, as the geodesics of $\Sph^3$ are great circles, the distance between $1:=(1,0,0,0)\in\Sph^3$ and a point $p:=(a,b,c,d)\in\Sph^3$ is given by $d_{\Sph^3}(1,p)=\arccos(a)$. Now, by \cite[Corollary 2.29]{lee_riema}, the induced distance $d_{\mathcal{F}_3(\R)}(1,\Pi(p))$ is obtained via
\[d_{g_8}(1,\Pi(p))=\min_{q\in\mathcal{Q}_8}d_{\Sph^3}(1,pq)=\min_{x=a,b,c,d}(\arccos|x|)=\arccos\|p\|_\infty\]
and therefore the diameter is
\[\diam(\mathcal{F}_3(\R),g_8)=\max_{\mathcal{F}_3(\R)}d_{g_8}(1,-)=\arccos\left(\min_{p\in\Sph^3}\|p\|_\infty\right)=\arccos(1/2)=\frac{\pi}{3}.\]
On the other hand, as already noted in Remark \ref{BWAC}, the diameter of the Killing metric on $\mathcal{F}_n$ has been computed in \cite{agrachev-chambrion} and is realized by the two 3-cycles of the symmetric groups. Because $\mathcal{F}_3(\R)\subset\mathcal{F}_3$ is a totally geodesic submanifold, we have
\[\diam(\mathcal{F}_3(\R),g_\kappa)=\diam(\mathcal{F}_3,g_\kappa)=\frac{\pi}{3}\sqrt{48}\]
and thus obtain
\begin{equation}\label{propto}
g_\kappa=48g_8.
\end{equation}

\begin{rem}\label{why}
In \cite{chirivi-garnier-spreafico}, the cells of the decomposition of $\mathcal{F}_3(\R)$ are obtained by gluing projections of cells of an $\mathcal{O}$-equivariant CW structure on $\Sph^3$. The latter cells are defined as \textit{curved joins} and thus are unions of open minimizing geodesics.

This means that the cells in $\mathcal{F}_3(\R)$ are unions of open geodesics for $g_\kappa$. In particular and recalling the notation from \S \ref{setting}, the $\Sym_3$-orbits of 1-cells are represented by the three cells
\[\mathfrak{e}_\delta:=\left\{\overline{\exp(tu_\delta)},~0<t<\pi/2\right\},\]
where $\delta\in\{\alpha,\beta,\alpha+\beta\}$ runs through the positive roots of the root system $\Phi=A_2$ (the upper bar denotes the image in $\mathcal{F}_3(\R)$ of an element of $SO(3)$). Observe that when $\delta\in\{\alpha,\beta\}$ is simple, the set
\[\left\{\overline{\exp(tu_\delta)},~0<t<\pi\right\}=\mathfrak{e}_\delta\sqcup\{s_\delta\}\sqcup(\mathfrak{e}_\delta \cdot s_\delta)\]
is the (real) \textit{Bruhat cell} associated to the simple reflection $s_\delta\in W$. More generally, the reader may observe that the real Bruhat cells are unions of equivariant cells in $\mathcal{F}_3(\R)$.

Finally, it should be noticed that the identity element is a vertex of the fundamental domain for $\Sym_3$ acting on $\mathcal{F}_3(\R)$ provided in \cite{chirivi-garnier-spreafico}, rather than its center. Moreover, the diameter of the domain exceeds twice the injectivity radius. 
\end{rem}

\subsection{Radius and polytopality of the Dirichlet--Voronoi domain in $\mathcal{F}_{A_2}(\R)$}

In this section, with the help of the quaternion algebra, we prove that the Dirichlet--Voronoi domain for $\mathcal{F}_3(\R)$ satisfies the hypotheses of Proposition \ref{cellsfromDV}, therefore providing an $\Sym_3$-equivariant cell structure. We thus view the sphere $\Sph^3$ as the space of quaternions with unit norm and we use the metric induced by $\Sph^3$ on $\mathcal{F}_3(\R)$ to do the calculations, as the associated distance is easier to handle. Throughout this section, we denote by
\[\mathcal{DV}_3=\{x\in\mathcal{F}_3(\R)~;~d(1,x)\le d(w,x),~\forall w\in\Sym_3\}\]
the Dirichlet--Voronoi domain for $\Sym_3$ acting on $\mathcal{F}_3(\R)$. We first determine the maximal value of the function $d(1,-)$ on $\mathcal{DV}_3$. Recall the Riemannian submersion $\Pi : \Sph^3\longtwoheadrightarrow\mathcal{F}_3(\R)$ from \eqref{projection_map} and that we denote by $s_\alpha=(1,2)$ and $s_\beta=(2,3)$ the generators of $W=\Sym_3$, together with the third reflection $s_{\alpha+\beta}:=s_\alpha s_\beta s_\alpha=s_\beta s_\alpha s_\beta$.

\begin{lem}\label{tcube}
Let $q=a+bi+cj+dk\in\Sph^3$ be such that $a\ge|b|,|c|,|d|$. Then
\[\Pi(q)\in\mathcal{DV}_3~\Longleftrightarrow~\left\{\begin{array}{rcl}|b|,|c|,|d| & \le & a(\sqrt{2}-1), \\[.5em] |b\pm c\pm d| & \le & a.\end{array}\right.\]
Moreover, letting $s_\gamma:=s_{\alpha+\beta}=s_\alpha s_\beta s_\alpha=s_\beta s_\alpha s_\beta$, we have
\begin{equation}\label{partial}\tag{$\partial$}
\left\{\begin{array}{rcl}
\Pi(q)\in Z_{s_\alpha} & \Longleftrightarrow & a(\sqrt{2}-1)=|d|, \\[.5em]
\Pi(q)\in Z_{s_\beta} & \Longleftrightarrow & a(\sqrt{2}-1)=|b|, \\[.5em]
\Pi(q)\in Z_{s_{\alpha+\beta}} & \Longleftrightarrow & a(\sqrt{2}-1)=|c|, \\[.5em]
\Pi(q)\in Z_{s_\alpha s_\beta} & \Longleftrightarrow & a=\max(b-c-d,-b+c-d,-b-c+d,b+c+d), \\[.5em]
\Pi(q)\in Z_{s_\beta s_\alpha} & \Longleftrightarrow & a=\max(b+c-d,b-c+d,-b+c+d,-b-c-d).\end{array}\right.
\end{equation}
\end{lem}
\begin{proof}
The hypothesis $a\ge|b|,|c|,|d|$ ensures that
\[d(1,\Pi(q))=4\sqrt{3}\arccos\|q\|_\infty=4\sqrt{3}\arccos(a).\]
As the elements $s_\alpha$, $s_\beta$, $s_{\alpha+\beta}=s_\alpha s_\beta s_\alpha$, $s_\alpha s_\beta$ and $s_\beta s_\alpha$ of $\Sym_3$ may be represented in $\Sph^3$ respectively by $\tfrac{1}{\sqrt{2}}(1+k)$, $\tfrac{1}{\sqrt{2}}(1+i)$, $\tfrac{1}{\sqrt{2}}(1+j)$, $\tfrac12(1+i+j+k)$ and $\tfrac12(1-i+j+k)$, we find
\[\left\{\begin{array}{rcl}
d(1,\Pi(q)s_\alpha) & = & 4\sqrt{3}\arccos\left(\tfrac{1}{\sqrt{2}}\max(|a\pm d|,|b\pm c|)\right), \\[.5em]
d(1,\Pi(q)s_\beta) & = & 4\sqrt{3}\arccos\left(\tfrac{1}{\sqrt{2}}\max(|a\pm b|,|c\pm d|)\right), \\[.5em]
d(1,\Pi(q)s_{\alpha+\beta}) & = & 4\sqrt{3}\arccos\left(\tfrac{1}{\sqrt{2}}\max(|a\pm c|,|b\pm d|)\right), \\[.5em]
d(1,\Pi(q)s_\alpha s_\beta) & = & 4\sqrt{3}\arccos\left(\tfrac12\max_{\epsilon_i=\pm1,~\epsilon_1\epsilon_2\epsilon_3=1}(|a+\epsilon_1b+\epsilon_2c+\epsilon_3d|)\right), \\[.5em]
d(1,\Pi(q)s_\beta s_\alpha) & = & 4\sqrt{3}\arccos\left(\tfrac12\max_{\epsilon_i=\pm1,~\epsilon_1\epsilon_2\epsilon_3=-1}(|a+\epsilon_1b+\epsilon_2c+\epsilon_3d|)\right).\end{array}\right.\]
Thus, the conditions $d(1,\Pi(q))\le d(w,\Pi(q))$ for $w\in\Sym_3\setminus\{1\}$ translate into the system of inequalities
\begin{equation}\label{sy}\tag{$S$}
\left\{\begin{array}{rcl}
a\sqrt{2} & \ge & \max(|a\pm d|,|b\pm c|,|a\pm b|,|c\pm d|,|a\pm c|,|b\pm d|), \\[.5em]
2a & \ge & \max(|a\pm b\pm c\pm d|).\end{array}\right.
\end{equation}
The second set of inequalities in \eqref{sy} is readily equivalent to the inequalities $|b\pm c\pm d|\le a$. On the other hand, since $|b|\le a$, we have $a\pm b\ge 0$ and therefore, the first set of inequalities in \eqref{sy} implies in particular that $a\sqrt{2}\ge a\pm b$, so that $a(\sqrt{2}-1)\ge |b|$. In the same way, $a(\sqrt{2}-1)\ge |c|,|d|$. The triangular inequality then yields $|b\pm c|\le|b|+|c|\le 2a(\sqrt{2}-1)<a\sqrt{2}$ and similarly, $|c\pm d|,|b\pm d|<a\sqrt{2}$. Therefore, the system \eqref{sy} is indeed equivalent to the system of the statement.

Now, using the formulae for the distance $d(1,\Pi(q)w)$ given above, we arrive at
\[\Pi(q)\in H_{s_\alpha}~\Longleftrightarrow~a\sqrt{2}=\max(|a\pm d|,|b\pm c|)=a+|d|~\Longleftrightarrow~a(\sqrt{2}-1)=|d|,\]
the conditions for $H_{s_\beta}$ and $H_{s_{\alpha+\beta}}$ being analogous. 
Observe that since $a\pm b\ge0$, we have the inequality $a-b-c+d\ge -a+b-c+d$ and similarly, $a+b-c-d\ge -a-b-c-d$, $a-b+c-d\ge -a+b+c-d$ and $a+b+c+d\ge -a-b+c+d$, so that
\begin{align*}
&\max(|a+b-c-d|,|a-b+c-d|,|a-b-c+d|,|a+b+c+d|) \\
=&\max(a+b-c-d,a-b+c-d,a-b-c+d,a+b+c+d) \\
=&a+\max(b-c-d,-b+c-d,-b-c+d,b+c+d).
\end{align*}
Thus, using the formula for $d(1,\Pi(q)s_\alpha s_\beta)$, we have equivalences
\begin{align*}
\Pi(q)\in H_{s_\alpha s_\beta}&\Longleftrightarrow 2a=\max(|a+b-c-d|,|a-b+c-d|,|a-b-c+d|,|a+b+c+d|) \\ 
&\Longleftrightarrow a=\max(b-c-d,-b+c-d,-b-c+d,b+c+d).
\end{align*}
\end{proof}

\begin{prop}\label{DVisokforSO(3)}
The radius
\[\max_{p\in\mathcal{DV}_3}d(1,p)=4\sqrt{3}\arccos\left(\frac{1}{2}+\frac{\sqrt{2}}{4}\right)=:\delta_0\approx3.7969\]
of $\mathcal{DV}_3$ is smaller than the injectivity radius $\pi\sqrt{3}\approx5.4414$ and in particular, the interior $\stackrel{\circ}{\mathcal{DV}_3}$ is a $3$-cell. Moreover, this maximum is attained by the twenty-four points
\[\Pi\left({\tfrac{1}{2}+\tfrac{\sqrt{2}}{4}+bi+cj+dk}\right),\]
where $(b,c,d)$ is any permutation of $\left(\pm\tfrac{\sqrt{2}}{4},\pm\tfrac{\sqrt{2}}{4},\pm\left(\tfrac{1}{2}-\tfrac{\sqrt{2}}{4}\right)\right)$.
\end{prop}

\begin{proof}
First, we have to prove that $\mathcal{DV}_3\subset B(1,\inj(\mathcal{F}_3(\R)))=B(1,\pi\sqrt{3})$. In view of Lemma \ref{tcube}, if $q=a+bi+cj+dk\in\mathcal{DV}_3$ (with $a\ge |b|,|c|,|d|$), then in particular $\max(|b|,|c|,|d|)\le a(\sqrt{2}-1)$, implying that
\[1-a^2=b^2+c^2+d^2\le 3a^2(\sqrt{2}-1)^2~~\Longrightarrow~~1\le a^2(10-6\sqrt{2})\]
and thus
\[d(1,\Pi({q}))=4\sqrt{3}\arccos(a)\le 4\sqrt{3}\arccos\left(\frac{1}{\sqrt{10-6\sqrt{2}}}\right)<4\sqrt{3}\arccos\left(\frac{1}{\sqrt{2}}\right)=\pi\sqrt{3},\]
as required.

Now that we know that $\mathcal{DV}_3\subset B(1,\inj(\mathcal{F}_3(\R)))$, in view of Lemma \ref{isolated_max}, it remains to show that there is at most a finite number of elements of $\mathcal{DV}_3$ at distance $\delta_0$ from $1$, and to find them. For $q=a+bi+cj+dk\in\mathcal{DV}_3$ as above, to say that $d(1,\Pi(q))=\delta_0$ amounts to say that $a=1/2+\sqrt{2}/4$. Therefore, the system of inequalities \eqref{sy} of the Lemma, together with the condition $a^2+b^2+c^2+d^2=1$, yields the system
\[\left\{\begin{array}{rcl}
|b|,|c|,|d| & \le & \sqrt{2}/4, \\[.5em]
|b\pm c\pm d| & \le & 1/2+\sqrt{2}/4, \\[.5em]
b^2+c^2+d^2 & = & \tfrac{1}{8}(5-2\sqrt{2}),\end{array}\right.\]
which defines a full truncated cube in $\R^3$, intersected with a sphere, as shown in Figure \ref{trunccube_red&blue}. 

\begin{center}
\begin{figure}[h!]
\begin{tikzpicture}[x={(0cm,-1cm)},y={(1cm,0cm)},z={(-2.85mm,-2.85mm)},scale=4.95]
	\tikzmath{function xx(\r,\t,\p) {return \r * sin(\t r) * sin(\p r);};}
	\tikzmath{function yy(\r,\t,\p) {return \r * sin(\t r) * cos(\p r);};}
	\tikzmath{function zz(\r,\t,\p) {return \r * cos(\t r);};}
	\pgfmathsetmacro{\radio}{sqrt(1/8*(5-2*sqrt(2)))}
	\pgfmathsetmacro{\Aa}{sqrt(2)/4}
	\pgfmathsetmacro{\bound}{\Aa/\radio}
	\pgfmathsetmacro{\boundb}{(1/2-\Aa)/\radio}
	
	\pgfmathsetmacro{\thetaa}{pi/180*acos(\bound)}
	\pgfmathsetmacro{\thetab}{pi/180*acos(\boundb)}
	
	\coordinate (ca) at (\Aa,\Aa,\Aa);
	\coordinate (cb) at (\Aa,\Aa,-\Aa);
	\coordinate (cc) at (\Aa,-\Aa,\Aa);
	\coordinate (cd) at (\Aa,-\Aa,-\Aa);
	\coordinate (ce) at (-\Aa,\Aa,\Aa);
	\coordinate (cf) at (-\Aa,\Aa,-\Aa);
	\coordinate (cg) at (-\Aa,-\Aa,\Aa);
	\coordinate (ch) at (-\Aa,-\Aa,-\Aa);
	
	\draw[dotted] (cd)--(ch) (cd)--(cb) (cd)--(cc);
	
	\coordinate (a) at (\Aa,\Aa,1/2-\Aa);
	\coordinate (b) at (\Aa,1/2-\Aa,\Aa);
	\coordinate (c) at (1/2-\Aa,\Aa,\Aa);
	\coordinate (d) at (\Aa,\Aa,-1/2+\Aa);
	\coordinate (e) at (\Aa,1/2-\Aa,-\Aa);
	\coordinate (f) at (1/2-\Aa,\Aa,-\Aa);
	\coordinate (g) at (\Aa,-\Aa,1/2-\Aa);
	\coordinate (h) at (\Aa,-1/2+\Aa,\Aa);
	\coordinate (i) at (1/2-\Aa,-\Aa,\Aa);
	\coordinate (j) at (-\Aa,\Aa,1/2-\Aa);
	\coordinate (k) at (-\Aa,1/2-\Aa,\Aa);
	\coordinate (l) at (-1/2+\Aa,\Aa,\Aa);
	
	\coordinate (m) at (-\Aa,-\Aa,-1/2+\Aa);
	\coordinate (n) at (-\Aa,-1/2+\Aa,-\Aa);
	\coordinate (o) at (-1/2+\Aa,-\Aa,-\Aa);
	\coordinate (p) at (-\Aa,-\Aa,1/2-\Aa);
	\coordinate (q) at (-\Aa,-1/2+\Aa,\Aa);
	\coordinate (r) at (-1/2+\Aa,-\Aa,\Aa);
	\coordinate (s) at (-\Aa,\Aa,-1/2+\Aa);
	\coordinate (t) at (-\Aa,1/2-\Aa,-\Aa);
	\coordinate (u) at (-1/2+\Aa,\Aa,-\Aa);
	\coordinate (v) at (\Aa,-\Aa,-1/2+\Aa);
	\coordinate (w) at (\Aa,-1/2+\Aa,-\Aa);
	\coordinate (x) at (1/2-\Aa,-\Aa,-\Aa);
	\fill[fill=blue,opacity=1] (g)--(h)--(i)--(g) (j)--(k)--(l)--(j) (p)--(q)--(r)--(p);
	
	\pgfmathsetmacro{\step}{0.005}
	\pgfmathsetmacro{\thastep}{\thetaa+\step}
	\pgfmathsetmacro{\thastepp}{\thastep+\step}
	\foreach \t in {\thetaa,\thastep,\thastepp,...,\thetab}{
		\pgfmathsetmacro{\eec}{pi/180*acos(\bound/sin(\t r))};
		\pgfmathsetmacro{\ees}{pi/180*asin(\bound/sin(\t r))};
		\draw[red,thin,opacity=1] plot[domain=\eec:\ees,smooth,variable=\p]
			({zz(\radio,\t,\p)},{yy(\radio,\t,\p)},{-xx(\radio,\t,\p)});
		\draw[red,thin,opacity=1] plot[domain=\eec:\ees,smooth,variable=\p]
			({zz(\radio,\t,\p)},{-yy(\radio,\t,\p)},{xx(\radio,\t,\p)});
		\draw[red,thin,opacity=1] plot[domain=\eec:\ees,smooth,variable=\p]
			({-zz(\radio,\t,\p)},{yy(\radio,\t,\p)},{-xx(\radio,\t,\p)});
		\draw[red,thin,opacity=1] plot[domain=\eec:\ees,smooth,variable=\p]
			({-zz(\radio,\t,\p)},{-yy(\radio,\t,\p)},{-xx(\radio,\t,\p)});
	}
	
	\fill[fill=blue,opacity=1] (a)--(b)--(c)--(a) (d)--(e)--(f)--(d) (m)--(n)--(o)--(m) (s)--(t)--(u)--(s) (v)--(w)--(x)--(v);
	
	\foreach \t in {\thetaa,\thastep,\thastepp,...,\thetab}{
		\pgfmathsetmacro{\eec}{pi/180*acos(\bound/sin(\t r))};
		\pgfmathsetmacro{\ees}{pi/180*asin(\bound/sin(\t r))};
		\draw[red,thin,opacity=1] plot[domain=\eec:\ees,smooth,variable=\p]
			({zz(\radio,\t,\p)},{yy(\radio,\t,\p)},{xx(\radio,\t,\p)});
		\draw[red,thin,opacity=1] plot[domain=\eec:\ees,smooth,variable=\p]
			({-zz(\radio,\t,\p)},{yy(\radio,\t,\p)},{xx(\radio,\t,\p)});
		\draw[red,thin,opacity=1] plot[domain=\eec:\ees,smooth,variable=\p]
			({-zz(\radio,\t,\p)},{-yy(\radio,\t,\p)},{xx(\radio,\t,\p)});
	}
	
	\draw (ca)--(cb) (ca)--(cc) (ca)--(ce) (ce)--(cf) (ce)--(cg) (cg)--(cc) (cb)--(cf) (cf)--(ch) (ch)--(cg);
\end{tikzpicture}
\caption{The blue triangles are faces of the truncated cube and the red parts form the intersection of the sphere with the full cube.}
\label{trunccube_red&blue}
\end{figure}
\end{center}

This full truncated (see Figure \ref{dirichlet_domain_R3}) cube is the convex hull of its twenty-four vertices we listed in the statement, and these lie on the sphere. As the Euclidean norm on $\R^3$ is uniformly convex, any point in this truncated cube which is not a vertex belongs to the open ball of radius $\sqrt{\tfrac{5-2\sqrt{2}}{8}}$, hence it is not on the sphere. Therefore, the considered intersection consists exactly of the twenty-four points of the statement.
\end{proof}

The previous result and proof suggest some kind of identification between $\mathcal{DV}_3$ and a truncated cube in $\R^3$. We make this more precise in the following result:
\begin{prop}\label{ident_tcube}
Let $\mathcal{K}\subset\R^3$ be the truncated cube whose vertices are all the permutations of the point $\left(\pm\tfrac{\sqrt{2}}{4},\pm\tfrac{\sqrt{2}}{4},\pm\left(\tfrac12-\tfrac{\sqrt{2}}{4}\right)\right)$ or equivalently, defined by the inequalities
\begin{equation}\label{eqK}\tag{$S_{\mathcal{K}}$}
\left\{\begin{array}{rcl}
|x|,|y|,|z| & \le & \sqrt{2}/4, \\[.5em]
|x\pm y\pm z| & \le & \sqrt{2}/4+1/2.\end{array}\right.
\end{equation}
Then, the map
\[\begin{array}{ccccc}
\varphi_{\mathcal{K}} & : & \mathcal{K} & \longto & \mathcal{DV}_3 \\ & & (x,y,z) & \longmapsto & \Pi\left(\frac{1+\sqrt{2}+2\sqrt{2}(xi+yj+zk)}{\sqrt{3+2\sqrt{2}+8x^2+8y^2+8z^2}}\right)\end{array}\]
defines a homeomorphism that sends the vertices of $\mathcal{K}$ to the points at maximal distance from $1$ in $\mathcal{DV}_3$, that is, to $\mathcal{DV}_3\cap S(1,\delta_0)$.
\end{prop}

\begin{center}
\begin{figure}[h!]
\begin{tikzpicture}[x={(0cm,-1cm)},y={(1cm,0cm)},z={(-2.85mm,-2.85mm)},scale=4.95]
	\pgfmathsetmacro{\radio}{sqrt(1/8*(5-2*sqrt(2)))}
	\pgfmathsetmacro{\Aa}{sqrt(2)/4}
	\pgfmathsetmacro{\bound}{\Aa/\radio}
	\pgfmathsetmacro{\boundb}{(1/2-\Aa)/\radio}
	
	\pgfmathsetmacro{\thetaa}{pi/180*acos(\bound)}
	\pgfmathsetmacro{\thetab}{pi/180*acos(\boundb)}
	
	\coordinate (ca) at (\Aa,\Aa,\Aa);
	\coordinate (cb) at (\Aa,\Aa,-\Aa);
	\coordinate (cc) at (\Aa,-\Aa,\Aa);
	\coordinate (cd) at (\Aa,-\Aa,-\Aa);
	\coordinate (ce) at (-\Aa,\Aa,\Aa);
	\coordinate (cf) at (-\Aa,\Aa,-\Aa);
	\coordinate (cg) at (-\Aa,-\Aa,\Aa);
	\coordinate (ch) at (-\Aa,-\Aa,-\Aa);
	
	\draw[dotted] (cd)--(ch) (cd)--(cb) (cd)--(cc);
	
	\coordinate (a) at (\Aa,\Aa,1/2-\Aa);
	\coordinate (b) at (\Aa,1/2-\Aa,\Aa);
	\coordinate (c) at (1/2-\Aa,\Aa,\Aa);
	\coordinate (d) at (\Aa,\Aa,-1/2+\Aa);
	\coordinate (e) at (\Aa,1/2-\Aa,-\Aa);
	\coordinate (f) at (1/2-\Aa,\Aa,-\Aa);
	\coordinate (g) at (\Aa,-\Aa,1/2-\Aa);
	\coordinate (h) at (\Aa,-1/2+\Aa,\Aa);
	\coordinate (i) at (1/2-\Aa,-\Aa,\Aa);
	\coordinate (j) at (-\Aa,\Aa,1/2-\Aa);
	\coordinate (k) at (-\Aa,1/2-\Aa,\Aa);
	\coordinate (l) at (-1/2+\Aa,\Aa,\Aa);
	
	\coordinate (m) at (-\Aa,-\Aa,-1/2+\Aa);
	\coordinate (n) at (-\Aa,-1/2+\Aa,-\Aa);
	\coordinate (o) at (-1/2+\Aa,-\Aa,-\Aa);
	\coordinate (p) at (-\Aa,-\Aa,1/2-\Aa);
	\coordinate (q) at (-\Aa,-1/2+\Aa,\Aa);
	\coordinate (r) at (-1/2+\Aa,-\Aa,\Aa);
	\coordinate (s) at (-\Aa,\Aa,-1/2+\Aa);
	\coordinate (t) at (-\Aa,1/2-\Aa,-\Aa);
	\coordinate (u) at (-1/2+\Aa,\Aa,-\Aa);
	\coordinate (v) at (\Aa,-\Aa,-1/2+\Aa);
	\coordinate (w) at (\Aa,-1/2+\Aa,-\Aa);
	\coordinate (x) at (1/2-\Aa,-\Aa,-\Aa);
	\fill[fill=pur,opacity=0.8] (a)--(b)--(c)--(a) (j)--(k)--(l)--(j);
	\fill[fill=greeo,opacity=0.8] (g)--(h)--(i)--(g) (p)--(q)--(r)--(p);		
	\fill[fill=greeo,opacity=0.8] (a)--(c)--(l)--(j)--(s)--(u)--(f)--(d)--(a);
	\fill[fill=cof,opacity=0.8] (b)--(c)--(l)--(k)--(q)--(r)--(i)--(h)--(b);
	\fill[fill=greet,opacity=0.8] (k)--(j)--(s)--(t)--(n)--(m)--(p)--(q)--(k);
	
	\draw[thick] (a)--(c) (c)--(l) (l)--(j) (j)--(s) (s)--(u) (u)--(f) (f)--(d) (d)--(a) (b)--(c) (c)--(l) (l)--(k) (k)--(q) (q)--(r) (r)--(i) (i)--(h) (h)--(b) (k)--(j) (j)--(s) (s)--(t) (t)--(n) (n)--(m) (m)--(p) (p)--(q) (q)--(k) (a)--(b) (r)--(p) (h)--(g) (g)--(i);
	
	\draw[thick, dashed] (g)--(v) (v)--(x) (v)--(w) (x)--(w) (x)--(o) (o)--(n) (o)--(m) (w)--(e) (e)--(f) (e)--(d) (t)--(u);
	
	\draw (ca)--(cb) (ca)--(cc) (ca)--(ce) (ce)--(cf) (ce)--(cg) (cg)--(cc) (cb)--(cf) (cf)--(ch) (ch)--(cg);
\end{tikzpicture}
\caption{The truncated cube $\mathcal{K}\simeq\mathcal{DV}_3$.}
\label{dirichlet_domain_R3}
\end{figure}
\end{center}

\begin{proof}
The Lemma \ref{tcube} ensures that $\varphi_{\mathcal{K}}$ is well-defined and the continuity is obvious. To prove that it is a homeomorphism, by compactness of $\mathcal{K}$, it suffices to show that $\varphi_{\mathcal{K}}$ is bijective.

Let $(x,y,z),(x',y',z')\in\mathcal{K}$ be such that $\varphi_{\mathcal{K}}(x,y,z)=\varphi_{\mathcal{K}}(x',y',z')$ and let $q,q'\in\Sph^3$ be the two quaternions considered in the formula defining $\varphi_{\mathcal{K}}$, i.e. such that $\varphi_{\mathcal{K}}(x,y,z)=\Pi(q)$ and $\varphi_{\mathcal{K}}(x',y',z')=\Pi(q')$. We have seen in the proof of Proposition \ref{DVisokforSO(3)} that $\mathcal{K}$ is included in the ball $B\left(0,\sqrt{\tfrac{5-2\sqrt{2}}{8}}\right)$, so that
\[\frac{1+\sqrt{2}}{\sqrt{3+2\sqrt{2}+8x^2+8y^2+8z^2}}\ge\frac{1+\sqrt{2}}{\sqrt{3+2\sqrt{2}+5-2\sqrt{2}}}=\frac12+\frac{\sqrt{2}}{4}.\]
In other words, the real parts of $q$ and $q'$ are greater than $1/2+\sqrt{2}/4$ and in particular, there is no $h\in\mathcal{Q}_8\setminus\{1\}=\{-1,\pm i,\pm j,\pm k\}$ such that $q'=qh$. Indeed, if for instance $q'=qi$, then writing $q=a+bi+cj+dk$ and $q'=qi=-b+ai+dj-ck$, we find that $|a|,|b|>1/2+\sqrt{2}/4$ and thus $1=a^2+b^2+c^2+d^2\ge a^2+b^2>2(1/2+\sqrt{2}/4)^2=3/4+\sqrt{2}/2>1$, a contradiction. The other cases are ruled out similarly and therefore, the equality $\Pi(q)=\Pi(q')$ implies $q=q'$. The real parts of $q$ and $q'$ being equal, we get $x^2+y^2+z^2=(x')^2+(y')^2+(z')^2$ and thus $(x,y,z)=(x',y',z')$, so that $\varphi_{\mathcal{K}}$ is one-to-one.

For the surjectivity, let $p\in\mathcal{DV}_3$ and choose coordinates $a,b,c,d$ such that $a\ge|b|,|c|,|d|$ and $p=\Pi(a+bi+cj+dk)$. By Lemma \ref{tcube}, the point $\frac{1/2+\sqrt{2}/4}{a}(b,c,d)$ lies in $\mathcal{K}$ and its image under $\varphi_{\mathcal{K}}$ is the image under $\Pi$ of the quaternion:
\begin{align*}
\frac{1+\sqrt{2}+\frac{1+\sqrt{2}}{a}(bi+cj+dk)}{\sqrt{3+2\sqrt{2}+\tfrac{8}{a^2}(1/2+\sqrt{2}/4)^2(b^2+c^2+d^2)}}&=\frac{(1+\sqrt{2})(a+bi+cj+dk)}{\sqrt{a^2(3+2\sqrt{2})+(3+2\sqrt{2})(1-a^2)}} \\
&=\frac{(1+\sqrt{2})(a+bi+cj+dk)}{\sqrt{3+2\sqrt{2}}} \\
&=a+bi+cj+dk,
\end{align*}
so that $\varphi_{\mathcal{K}}\left(\tfrac{1/2+\sqrt{2}/4}{a}(b,c,d)\right)=\Pi(a+bi+cj+dk)=p$, and thus $\varphi_{\mathcal{K}}$ is onto as well. The fact that the vertices of $\mathcal{K}$ are sent to $\mathcal{DV}_3\cap S(1,\delta_0)$ is immediate from the second statement of Proposition \ref{DVisokforSO(3)}.
\end{proof}

\begin{cor}\label{ident_faces}
For symbols $\epsilon_1,\epsilon_2,\epsilon_3\in\{+,-\}$ and $t\in\{x,y,z\}$, we denote by $F^{\epsilon_1,\epsilon_2,\epsilon_3}$ and $F^{\epsilon_1}_t$ the facets of $\mathcal{K}$ defined by
\[F^{\epsilon_1,\epsilon_2,\epsilon_3}:=\{(x,y,z)\in\mathcal{K}~|~\epsilon_1x+\epsilon_2y+\epsilon_3z=\sqrt{2}/4+1/2\},~F^{\epsilon_1}_t:=\{(x,y,z)\in\mathcal{K}~|~\epsilon_1t=\sqrt{2}/{4}\}.\]
Then, we have
\[\left\{\begin{array}{rcl}
\varphi_{\mathcal{K}}^{-1}(Z_{s_\alpha}) & = & F_z^+\sqcup F_z^-, \\[.5em]
\varphi_{\mathcal{K}}^{-1}(Z_{s_\beta}) & = & F_x^+\sqcup F_x^-, \\[.5em]
\varphi_{\mathcal{K}}^{-1}(Z_{s_{\alpha+\beta}}) & = & F_y^+\sqcup F_y^-, \\[.5em]
\varphi_{\mathcal{K}}^{-1}(Z_{s_\alpha s_\beta}) & = & F^{+,-,-}\sqcup F^{-,+,-}\sqcup F^{-,-,+}\sqcup F^{+,+,+}=\bigcup_{\epsilon_1\epsilon_2\epsilon_3=+}F^{\epsilon_1,\epsilon_2,\epsilon_3}, \\[.5em]
\varphi_{\mathcal{K}}^{-1}(Z_{s_\beta s_\alpha}) & = & F^{+,+,-}\sqcup F^{+,-,+}\sqcup F^{-,+,+}\sqcup F^{-,-,-}=\bigcup_{\epsilon_1\epsilon_2\epsilon_3=-}F^{\epsilon_1,\epsilon_2,\epsilon_3}.\end{array}\right.\]
\end{cor}
\begin{proof}
This follows easily from Lemma \ref{tcube} and the definition of $\varphi_{\mathcal{K}}$, keeping in mind that the system \eqref{partial} is invariant under positive dilatations.
\end{proof}

\subsection{The cell structure and its cellular chain complex}

We are now in the position to apply the general results of the first part to the case of $\mathcal{F}_3(\R)$. First, Propositions \ref{DVisokforSO(3)}, \ref{ident_tcube} and Corollary \ref{ident_faces} allow us to apply Proposition \ref{cellsfromDV} to obtain the following result:
\begin{cor}\label{indeedequivdec}
Denoting by $F(\mathcal{K})$ the face lattice of the truncated cube $\mathcal{K}$, the decomposition
\[\mathcal{F}_3(\R)=\coprod_{\substack{f\in F(\mathcal{K}) \\ w\in\Sym_3}}\varphi_{\mathcal{K}}(\interior{f})w\]
is an $\Sym_3$-equivariant CW structure on $\mathcal{F}_3(\R)$.
\end{cor}

To derive the associated cellular chain complex, we have to find representatives for the orbits of the cells and to compute their boundary. For this, we make heavy use of Lemma \ref{partactonwalls}. For instance, if $\varphi_{\mathcal{K}}(\interior{F}^{+,+,+})w\cap\mathcal{DV}_3\ne\emptyset$ for some $1\ne w\in\Sym_3$, then $w=s_\beta s_\alpha$ and $\varphi_{\mathcal{K}}(\interior{F}^{+,+,+})s_\beta s_\alpha=\interior{F}^{\epsilon_1,\epsilon_2,\epsilon_3}$ for some sign triplet $(\epsilon_1,\epsilon_2,\epsilon_3)$ satisfying $\epsilon_1\epsilon_2\epsilon_3=-1$. Since $\interior{F}^{+,+,+}$ is the relative interior of the convex hull of its bounding vertices $\{v_j\}_j$, we only have to look for the facet $F^{\epsilon_1,\epsilon_2,\epsilon_3}\in F(\mathcal{K})$ having $\{\varphi_{\mathcal{K}}^{-1}(\varphi_{\mathcal{K}}(v_j)s_\beta s_\alpha)\}_j$ as vertices. We find that $\varphi_{\mathcal{K}}(\interior{F}^{+,+,+})s_\beta s_\alpha=\interior{F}^{-,-,-}$, and this generalizes to
\[\varphi_{\mathcal{K}}(\interior{F}^{\epsilon_1,\epsilon_2,\epsilon_3})s_\beta s_\alpha=\interior{F}^{-\epsilon_1,-\epsilon_2,-\epsilon_3}\]
for each sign triplet $(\epsilon_i)$ such that $\prod_i\epsilon_i=+1$. Also, if $s\in \Sym_3$ is a reflection, then $Z_s=\varphi_{\mathcal{K}}(F^+_t)\sqcup\varphi_{\mathcal{K}}(F^-_t)$ for some $t\in\{x,y,z\}$ and we have
\[\varphi_{\mathcal{K}}(F^{\pm}_t)s=\varphi_{\mathcal{K}}(F^{\mp}_t).\]
Hence, there are seven orbits of 2-cells.

We do the same for the 1-cells with the $Z_u\cap Z_v$'s, noticing that $Z_{s_\alpha s_\beta}\cap Z_{s_\beta s_\alpha}=\emptyset$. All the other intersections are non-empty and decompose as a union of four closed 1-cells. Moreover,
\[(Z_{s_\alpha}\cap Z_{s_{\alpha+\beta}})s_\alpha s_\beta s_\alpha=Z_{s_{\alpha+\beta}}\cap Z_{s_\beta s_\alpha}\text{ and }(Z_{s_\beta}\cap Z_{s_{\alpha+\beta}})s_{\alpha}s_\beta s_\alpha=Z_{s_{\alpha+\beta}}\cap Z_{s_\alpha s_\beta},\]
as well as
\[(Z_{s_\beta}\cap Z_{s_{\alpha+\beta}})s_\beta= Z_{s_\beta}\cap Z_{s_\beta s_\alpha}\text{ and }(Z_{s_\alpha}\cap Z_{s_{\beta}})s_\beta=Z_{s_\beta}\cap Z_{s_\alpha s_\beta},\]
and also
\[(Z_{s_\alpha}\cap Z_{s_\beta})s_\alpha=Z_{s_\alpha}\cap Z_{s_\beta s_\alpha}\text{ and }(Z_{s_\alpha}\cap Z_{s_{\alpha+\beta}})s_\alpha=Z_{s_\alpha}\cap Z_{s_\alpha s_\beta}.\]
Thus, the representatives of orbits of the 1-cells may be chosen among the three intersections $Z_{s_\alpha}\cap Z_{s_\beta}$, $Z_{s_\alpha}\cap Z_{s_{\alpha+\beta}}$ and $Z_{s_\beta}\cap Z_{s_{\alpha+\beta}}$. This yields twelve orbits of 1-cells.

Then, we have $Z_{s_\alpha}\cap Z_{s_\beta}\cap Z_{s_{\alpha+\beta}}=\emptyset$ and the six non-empty intersections $Z_u\cap Z_v\cap Z_w$ each contain (the images of) four vertices of $\mathcal{K}$. Moreover, if $x\in Z_u\cap Z_v\cap Z_w$ is a vertex, then $x$, $xu^{-1}$, $xv^{-1}$ and $xw^{-1}$ are four vertices of $\mathcal{DV}$, belonging to four different intersections of three walls. There are thus six orbits of 0-cells.

We choose the representing cells ${e}^i_j:=\varphi_{\mathcal{K}}(\interior{f}^i_j)$, where $f^i_j$ are the following faces of $\mathcal{K}$:
\[\hspace{-.5cm}\left\{\begin{array}{rcl}
{f}^2_1 & = & F_x^+, \\[.5em]
{f}^2_2 & = & F_z^+, \\[.5em]
{f}^2_3 & = & F_y^+, \\[.5em]
{f}^2_4 & = & F^{+,-,+}, \\[.5em]
{f}^2_5 & = & F^{+,-,-}, \\[.5em]
{f}^2_6 & = & F^{+,+,-}, \\[.5em]
{f}^2_7 & = & F^{+,+,+}.\end{array}\right. \left\{\begin{array}{rcl}
{f}^1_1 & = & F_x^+\cap F_y^-, \\[.5em]
{f}^1_2 & = & F_x^+\cap F^{+,-,-}, \\[.5em]
{f}^1_3 & = & F_x^+\cap F_z^-, \\[.5em]
{f}^1_4 & = & F_x^+\cap F^{+,+,-}, \\[.5em]
{f}^1_5 & = & F_x^+\cap F_y^+, \\[.5em]
{f}^1_6 & = & F_x^+\cap F^{+,+,+}, \end{array}\right. \left\{\begin{array}{rcl}
{f}^1_7 & = & F_x^+\cap F_z^+, \\[.5em]
{f}^1_8 & = & F_x^+\cap F^{+,-,+}, \\[.5em]
{f}^1_9 & = & F_y^-\cap F_z^+, \\[.5em]
{f}^1_{10} & = & F_y^+\cap F_z^-, \\[.5em]
{f}^1_{11} & = & F_y^-\cap F^{+,-,+}, \\[.5em]
{f}^1_{12} & = & F_y^+\cap F^{+,+,-}.\end{array}\right.\] \[\left\{\begin{array}{rcl}
{f}^0_1 & = & F_y^+\cap F^{+,+,-}\cap F_z^-, \\[.5em]
{f}^0_2 & = & F_x^+\cap F^{+,+,-}\cap F_y^+, \\[.5em]
{f}^0_3 & = & F_y^+\cap F^{+,+,+}\cap F_z^+, \\[.5em]
{f}^0_4 & = & F_x^-\cap F^{-,+,+}\cap F_z^+, \\[.5em]
{f}^0_5 & = & F_y^+\cap F^{-,+,+}\cap F_z^+, \\[.5em]
{f}^0_6 & = & F_x^+\cap F^{+,+,-}\cap F_z^-.\end{array}\right.\]
Of course, the maximal cells are represented by $e^3:=\interior{\mathcal{DV}}=\varphi_{\mathcal{K}}(\interior{\mathcal{K}})$.
To compute the oriented boundaries, we choose to orient $\mathcal{K}\subset\R^3$ directly and the representing 2-cells accordingly; the resulting pictures are given in Appendix \ref{truncated_cube}. We finally arrive to the following main result:

\begin{theo}\label{last_decomposition}
The cellular homology chain complex associated to the $\Sym_3$-equivariant CW structure given in Corollary \ref{indeedequivdec} is isomorphic to the complex of $\Z[\Sym_3]$-modules
\[\xymatrix{\Z[\Sym_3] \ar^{\partial_3}[r] & \Z[\Sym_3]^7 \ar^{\partial_2}[r] & \Z[\Sym_3]^{12} \ar^{\partial_1}[r] & \Z[\Sym_3]^6}\]
whose boundaries are given by (left) multiplication by the matrices
\begin{align*}
&\partial_1=\left(\begin{smallmatrix}0 & 0 & 0 & 0 & 0 & s_\beta & -s_\beta & 0 & 0 & -1 & 0 & 1 \\ 0 & 0 & 0 & 1 & -1 & 0 & 0 & 0 & s_{\beta}s_\alpha & 0 & 0 & -1 \\ -s_\gamma & 0 & 0 & 0 & 0 & 0 & 0 & s_\gamma & 0 & s_\beta & -s_\gamma & 0 \\ s_\beta s_\alpha & -s_\beta s_\alpha & 0 & 0 & s_\alpha & -s_\alpha & 0 & 0 & 0 & 0 & 0 & 0 \\ 0 & s_\beta s_\alpha & -s_\beta s_\alpha & 0 & 0 & 0 & 0 & 0 & -s_\gamma & 0 & s_\gamma & 0 \\ 0 & 0 & 1 & -1 & 0 & 0 & s_\beta & -s_\beta & 0 & 0 & 0 & 0\end{smallmatrix}\right), \\
&\partial_2=\left(\begin{smallmatrix}
1 & 0 & s_\gamma & 0 & 0 & 0 & -s_\gamma \\
1 & -s_\alpha s_\beta & 0 & 0 & -1 & 0 & 0 \\
1 & s_\beta & 0 & -s_\beta & 0 & 0 & 0 \\
1 & 0 & s_\alpha & 0 & 0 & -1 & 0 \\
1 & 0 & -1 & 0 & -s_\gamma & 0 & 0 \\
1 & s_\alpha & 0 & 0 & 0 & 0 & -1 \\
1 & -1 & 0 & 0 & 0 & -s_\beta & 0 \\
1 & 0 & -s_\beta s_\alpha & -1 & 0 & 0 & 0 \\
0 & -1 & -s_\gamma & 0 & -s_\beta & 0 & 0 \\
0 & s_\beta & 1 & 0 & 0 & 0 & -s_\beta \\
0 & s_\gamma & -s_\gamma & -1 & 0 & 0 & 0 \\
0 & -s_\beta s_\alpha & 1 & 0 & 0 & -1 & 0\end{smallmatrix}\right),~
\partial_3=\left(\begin{smallmatrix}
1-s_\alpha \\
1-s_\beta \\
1-s_\gamma \\
1-s_\beta s_\alpha \\
1-s_\alpha s_\beta \\
1-s_\beta s_\alpha \\
1-s_\alpha s_\beta\end{smallmatrix}\right),
\end{align*}
where $s_\alpha=(1,2)$, $s_\beta=(2,3)$ and $s_\gamma:=s_{\alpha+\beta}=s_\alpha s_\beta s_\alpha=s_\beta s_\alpha s_\beta$.
\end{theo}

Above, we chose the representing cells so that the situation can be easily visualized in Appendix \ref{truncated_cube}. However, it may seem more natural to choose the triangular faces in a single wall; the wall $Z_{s_\alpha s_\beta}$ for instance. If we let $\widetilde{e}^i_j=\varphi_{\mathcal{K}}(\interior{\widetilde{f}}^i_j)$, where
\scriptsize{\[\hspace{-.5cm}\left\{\begin{array}{rcl}
\widetilde{f}^2_1 & = & F_z^+, \\[.5em]
\widetilde{f}^2_2 & = & F_x^+, \\[.5em]
\widetilde{f}^2_3 & = & F_y^+, \\[.5em]
\widetilde{f}^2_4 & = & F^{+,+,+}, \\[.5em]
\widetilde{f}^2_5 & = & F^{+,-,-}, \\[.5em]
\widetilde{f}^2_6 & = & F^{-,+,-}, \\[.5em]
\widetilde{f}^2_7 & = & F^{-,-,+}.\end{array}\right. \left\{\begin{array}{rcl}
\widetilde{f}^1_1 & = & F_y^+\cap F^{+,+,+}, \\[.5em]
\widetilde{f}^1_2 & = & F_x^+\cap F^{+,-,-}, \\[.5em]
\widetilde{f}^1_3 & = & F_x^+\cap F_z^-, \\[.5em]
\widetilde{f}^1_4 & = & F_x^-\cap F_y^+, \\[.5em]
\widetilde{f}^1_5 & = & F_x^+\cap F_y^+, \\[.5em]
\widetilde{f}^1_6 & = & F_x^+\cap F^{+,+,+}, \end{array}\right. \left\{\begin{array}{rcl}
\widetilde{f}^1_7 & = & F_x^+\cap F_z^+, \\[.5em]
\widetilde{f}^1_8 & = & F_y^+\cap F^{-,+,-}, \\[.5em]
\widetilde{f}^1_9 & = & F_y^-\cap F_z^+, \\[.5em]
\widetilde{f}^1_{10} & = & F_z^+\cap F^{+,+,+}, \\[.5em]
\widetilde{f}^1_{11} & = & F_y^+\cap F_z^+, \\[.5em]
\widetilde{f}^1_{12} & = & F_z^+\cap F^{-,-,+}.\end{array}\right. \left\{\begin{array}{rcl}
\widetilde{f}^0_1 & = & F_x^+\cap F^{+,+,+}\cap F_z^+, \\[.5em]
\widetilde{f}^0_2 & = & F_x^+\cap F^{+,+,-}\cap F_y^+, \\[.5em]
\widetilde{f}^0_3 & = & F_y^+\cap F^{+,+,+}\cap F_z^+, \\[.5em]
\widetilde{f}^0_4 & = & F_x^+\cap F^{+,+,+}\cap F_y^+, \\[.5em]
\widetilde{f}^0_5 & = & F_y^+\cap F^{-,+,+}\cap F_z^+, \\[.5em]
\widetilde{f}^0_6 & = & F_x^+\cap F^{+,-,+}\cap F_z^+,\end{array}\right.\]}
\normalsize{then} the matrices of the boundaries $\widetilde{\partial}_i$ become
\begin{align*}
&\widetilde{\partial}_1=\left(\begin{smallmatrix}0 & 0 & 0 & 0 & 0 & 1 & -1 & 0 & 0 & -1 & 0 & s_\beta \\
0 & 0 & 0 & s_\beta & -1 & 0 & 0 & 0 & s_\alpha s_\beta & 0 & 0 & -s_\alpha s_\beta \\
-1 & 0 & 0 & 0 & 0 & 0 & 0 & s_\alpha & 0 & 1 & -1 & 0 \\
1 & -s_\gamma & 0 & 0 & 1 & -1 & 0 & 0 & 0 & 0 & 0 & 0 \\
0 & s_\alpha s_\beta & -s_\alpha s_\beta & 0 & 0 & 0 & 0 & 0 & -s_\gamma & 0 & 1 & 0 \\
0 & 0 & s_\alpha & -s_\alpha s_\beta & 0 & 0 & 1 & -s_\alpha s_\beta & 0 & 0 & 0 & 0\end{smallmatrix}\right), \\[.5em]
&\widetilde{\partial}_2=\left(\begin{smallmatrix}
0 & s_\gamma & 1 & -1 & 0 & 0 & 0 \\
-s_\beta s_\alpha & 1 & 0 & 0 & -1 & 0 & 0 \\
s_\alpha & 1 & 0 & 0 & 0 & s_\beta & 0 \\
0 & s_\beta & 1 & 0 & 0 & 0 & s_\gamma \\
0 & 1 & -1 & 0 & -s_\gamma & 0 & 0 \\
s_\beta & 1 & 0 & -1 & 0 & 0 & 0 \\
-1 & 1 & 0 & 0 & 0 & 0 & s_\beta \\
0 & s_\beta s_\alpha & -1 & 0 & 0 & 1 & 0 \\
-1 & 0 & -s_\gamma & 0 & -s_\alpha & 0 & 0 \\
1 & 0 & s_\alpha & -1 & 0 & 0 & 0 \\
1 & 0 & -1 & 0 & 0 & s_\alpha & 0 \\
-1 & 0 & s_\beta s_\alpha & 0 & 0 & 0 & 1\end{smallmatrix}\right),~~\widetilde{\partial}_3=\left(\begin{smallmatrix}
1-s_\alpha \\
1-s_\beta \\
1-s_\gamma \\
1-s_\beta s_\alpha \\
1-s_\beta s_\alpha \\
1-s_\beta s_\alpha \\
1-s_\beta s_\alpha\end{smallmatrix}\right).
\end{align*}

\begin{rem}
The homeomorphism $\mathcal{K}\stackrel{\tiny{\sim}}\to\mathcal{DV}_3$ crucially depends on the quaternion covering $\Sph^3\twoheadrightarrow SO(3)$. However, the inverse image $\mathrm{Exp}^{-1}(\mathcal{DV}_3)\subset\mathfrak{so}(3)$ has a convex hull
\[\mathcal{K}_{\mathfrak{so}}:=\mathrm{conv}(\mathrm{Exp}^{-1}(\mathcal{DV}_3)),\]
which is combinatorially equivalent to $\mathcal{K}$. Moreover, each inverse $I$-wall $\mathrm{Exp}^{-1}(Z_I)$ projects onto a union of faces of $\mathcal{K}_{\mathfrak{so}}$, corresponding to the faces appearing in $\varphi_{\mathcal{K}}^{-1}(Z_I)$. Therefore, the combinatorial information carried by $\mathcal{K}$ can be seen intrinsically in the Lie algebra $\mathfrak{so}(3)$.

To see this, recall the covering $\mu : \Sph^3\twoheadrightarrow SO(3)$ and take a quaternion $q=a_0+bi+cj+dk$ such that $a_0:=\tfrac12+\tfrac{\sqrt{2}}{4}\ge |b|,|c|,|d|$ (so that $\Pi(q)\in S(1,\delta_0)$). We have $\tr(\mu(q))=4a^2-1=\tfrac12+\sqrt{2}$ and if 
\[\theta_0:=\arccos\left(\tfrac{\tr(\mu(q))-1}{2}\right)=\arccos\left(\tfrac{\sqrt{2}}{2}-\tfrac14\right)\in]\tfrac\pi4,\tfrac\pi2[,\]
then, using the \textit{Rodrigues formula} \cite[\S 2]{log_rotations},
\[\log(\mu(q))=\frac{\theta_0}{2\sin\theta_0}\left(\mu(q)-{}^t{\mu(q)}\right)=\frac{\theta_0}{2\sin\theta_0}\begin{pmatrix}0 & -4a_0d & 4a_0c \\ 4a_0d & 0 & -4a_0b \\ -4a_0c & 4a_0b & 0\end{pmatrix}=k_0\begin{pmatrix}0 & -d & c \\ d & 0 & -b \\ -c & b & 0\end{pmatrix},\]
where
\[k_0:=\frac{4a_0\theta_0}{2\sin\theta_0}=\frac{2(2+\sqrt{2})\arccos(\sqrt{2}/2-1/4)}{\sqrt{7+4\sqrt{2}}}>0.\]
Observe also that if $(x,y,z)\in\mathcal{K}$ is a vertex of $\mathcal{K}$, then $(x,y,z)$ belongs to the sphere $\mathcal{S}_0$ centered at the origin and with radius $\sqrt{\tfrac18(5-2\sqrt{2})}$. Thus $\varphi_{\mathcal{K}}(x,y,z)=\Pi(a_0+xi+yj+zk)$ and so
\[\mathrm{Exp}^{-1}(\varphi_{\mathcal{K}}(x,y,z))=\log(\mu(a_0+xi+yj+zk))=k_0\begin{pmatrix}0 & -z & y \\ z & 0 & -x \\ -y & x & 0\end{pmatrix}.\]
Hence, if $V$ denotes the set of vertices of $\mathcal{K}$, then the assignment
\[\begin{array}{ccc}\mathcal{K} & \longto & \mathrm{conv}(\mathrm{Exp}^{-1}(\varphi_{\mathcal{K}}(V)))\subset\mathfrak{so}(3) \\\sum_{p}\lambda_pp & \longmapsto & \sum_p\lambda_p\mathrm{Exp}^{-1}(\varphi_{\mathcal{K}}(p))\end{array}\]
is an isomorphism of polytopes, where $\sum_p\lambda_pp$ denotes any convex combination of the vertices $p\in V$. On the other hand, each radial segment joining $0\in\mathfrak{so}(3)$ to a point on the sphere $S_{\mathfrak{so}(3)}(0,\delta_0)$ intersects exactly one inverse cell $\mathrm{Exp}^{-1}(e)\subset\mathrm{Exp}^{-1}(Z_I)\subset \mathrm{Exp}^{-1}(\partial\mathcal{DV}_3)$ at exactly one point. Furthermore, this segment also intersects exactly one boundary face of $\mathrm{conv}(\mathrm{Exp}^{-1}(\varphi_{\mathcal{K}}(V)))$ at exactly one point. This leads to a second isomorphism of polytopes
\[\mathrm{conv}(\mathrm{Exp}^{-1}(\varphi_{\mathcal{K}}(V)))\stackrel{\tiny{\sim}}\longto\mathrm{conv}(\mathrm{Exp}^{-1}(\mathcal{DV}_3))\stackrel{\tiny{\text{df}}}=\mathcal{K}_{\mathfrak{so}},\]
finally yielding an isomorphism of polytopes
\[\mathcal{K}\stackrel{\tiny{\sim}}\longto\mathcal{K}_{\mathfrak{so}}.\]
\end{rem}

\begin{appendix}
\addcontentsline{toc}{part}{Appendices}
\section{Combinatorics of the Dirichlet--Voronoi domain for $\mathcal{F}_{A_2}(\R)$}\label{truncated_cube}
Here we give some figures (Figures \ref{fig7} and \ref{fig8}) that help visualize how to obtain the complex from Theorem \ref{last_decomposition}. The cells belonging to the same $\Sym_3$-orbit share the same color and the order on the colors in the legends corresponds to the order chosen to build the matrices of the boundaries of the complex. For simplicity we replace $s_\alpha$, $s_\beta$ and $s_\gamma=s_\alpha s_\beta s_\alpha$ respectively by $a$, $b$ and $c$.

\begin{figure}[h!]
\begin{subfigure}[h]{0.3\textwidth}
\centering
\begin{tikzpicture}[x={(0cm,-1cm)},y={(1cm,0cm)},z={(-3.85mm,-3.85mm)},scale=1.75]
	\coordinate (a) at (0.4142,1,1);
	\coordinate (b) at (0.4142,1,-1);
	\coordinate (c) at (0.4142,-1,1);
	\coordinate (d) at (0.4142,-1,-1);
	\coordinate (e) at (-0.4142,1,1);
	\coordinate (f) at (-0.4142,1,-1);
	\coordinate (g) at (-0.4142,-1,1);
	\coordinate (h) at (-0.4142,-1,-1);
	\coordinate (i) at (1,0.4142,1);
	\coordinate (j) at (-1,0.4142,1);
	\coordinate (k) at (1,0.4142,-1);
	\coordinate (l) at (-1,0.4142,-1);
	\coordinate (m) at (1,-0.4142,1);
	\coordinate (n) at (-1,-0.4142,1);
	\coordinate (o) at (1,-0.4142,-1);
	\coordinate (p) at (-1,-0.4142,-1);
	\coordinate (q) at (1,1,0.4142);
	\coordinate (r) at (1,-1,0.4142);
	\coordinate (s) at (-1,1,0.4142);
	\coordinate (t) at (-1,-1,0.4142);
	\coordinate (u) at (1,1,-0.4142);
	\coordinate (v) at (1,-1,-0.4142);
	\coordinate (w) at (-1,1,-0.4142);
	\coordinate (x) at (-1,-1,-0.4142);

	\draw (a)--(i)--(m)--(c)--(g)--(n)--(j)--(e)--(a);
	\draw (b)--(f);
	\draw[dashed,opacity=0.5] (f)--(l);
	\draw[dashed,opacity=0.5] (p)--(h)--(d)--(o);
	\draw[dashed,opacity=0.5] (k)--(b);
	\draw (a)--(e)--(s)--(w)--(f)--(b)--(u)--(q)--(a);
	\draw (c)--(g)--(t)--(x);
	\draw[dashed,opacity=0.5] (x)--(h)--(d)--(v)--(r);
	\draw (r)--(c);
	\draw (n)--(t)--(x)--(p)--(l)--(w)--(s)--(j);
	\draw (i)--(m)--(r);
	\draw[dashed,opacity=0.5] (r)--(v)--(o)--(k)--(u);
	\draw (u)--(q)--(i);

	\fill[fill=black] (q) circle (1.5pt);
	\fill[fill=black] (v) circle (1.5pt);
	\fill[fill=black] (p) circle (1.5pt);
	\fill[fill=black] (j) circle (1.5pt);
	
	\fill[fill=red] (a) circle (1.5pt);
	\fill[fill=red] (f) circle (1.5pt);
	\fill[fill=red] (r) circle (1.5pt);
	\fill[fill=red] (x) circle (1.5pt);
	
	\fill[fill=cyan] (s) circle (1.5pt);
	\fill[fill=cyan] (u) circle (1.5pt);
	\fill[fill=cyan] (d) circle (1.5pt);
	\fill[fill=cyan] (g) circle (1.5pt);
	
	\fill[fill=green] (l) circle (1.5pt);
	\fill[fill=green] (e) circle (1.5pt);
	\fill[fill=green] (c) circle (1.5pt);
	\fill[fill=green] (o) circle (1.5pt);
	
	\fill[fill=brown] (w) circle (1.5pt);
	\fill[fill=brown] (k) circle (1.5pt);
	\fill[fill=brown] (t) circle (1.5pt);
	\fill[fill=brown] (m) circle (1.5pt);
	
	\fill[fill=pink] (i) circle (1.5pt);
	\fill[fill=pink] (n) circle (1.5pt);
	\fill[fill=pink] (h) circle (1.5pt);
	\fill[fill=pink] (b) circle (1.5pt);

	\draw (a) node[left]{$1$};
	\draw (b) node[right]{$a$};
	\draw (c) node[left]{$ba$};
	\draw (d) node[left]{$ab$};
	\draw (e) node[right]{$a$};
	\draw (f) node[right]{$a$};
	\draw (g) node[left]{$c$};
	\draw (h) node[right]{$ba$};
	\draw (i) node[below]{$1$};
	\draw (j) node[above]{$b$};
	\draw (k) node[below]{$b$};
	\draw (l) node[above]{$1$};
	\draw (m) node[below]{$ba$};
	\draw (n) node[above]{$b$};
	\draw (o) node[below]{$b$};
	\draw (p) node[above]{$ba$};
	\draw (q) node[below right]{$1$};
	\draw (r) node[below left]{$c$};
	\draw (s) node[above]{$1$};
	\draw (t) node[left]{$c$};
	\draw (u) node[below right]{$b$};
	\draw (v) node[below right]{$c$};
	\draw (w) node[above right]{$1$};
	\draw (x) node[above left]{$ba$};
\end{tikzpicture}
\end{subfigure}
\hspace{3cm}
\begin{subfigure}[h]{0.3\textwidth}
\begin{tikzpicture}[x={(0cm,-1cm)},y={(1cm,0cm)},z={(-3.85mm,-3.85mm)},scale=1.75]
	\coordinate (a) at (0.4142,1,1);
	\coordinate (b) at (0.4142,1,-1);
	\coordinate (c) at (0.4142,-1,1);
	\coordinate (d) at (0.4142,-1,-1);
	\coordinate (e) at (-0.4142,1,1);
	\coordinate (f) at (-0.4142,1,-1);
	\coordinate (g) at (-0.4142,-1,1);
	\coordinate (h) at (-0.4142,-1,-1);
	\coordinate (i) at (1,0.4142,1);
	\coordinate (j) at (-1,0.4142,1);
	\coordinate (k) at (1,0.4142,-1);
	\coordinate (l) at (-1,0.4142,-1);
	\coordinate (m) at (1,-0.4142,1);
	\coordinate (n) at (-1,-0.4142,1);
	\coordinate (o) at (1,-0.4142,-1);
	\coordinate (p) at (-1,-0.4142,-1);
	\coordinate (q) at (1,1,0.4142);
	\coordinate (r) at (1,-1,0.4142);
	\coordinate (s) at (-1,1,0.4142);
	\coordinate (t) at (-1,-1,0.4142);
	\coordinate (u) at (1,1,-0.4142);
	\coordinate (v) at (1,-1,-0.4142);
	\coordinate (w) at (-1,1,-0.4142);
	\coordinate (x) at (-1,-1,-0.4142);
	
	\coordinate (ae) at ($(a)!0.5!(e)$);
	\coordinate (ej) at ($(e)!0.5!(j)$);
	\coordinate (jn) at ($(j)!0.5!(n)$);
	\coordinate (ng) at ($(n)!0.5!(g)$);
	\coordinate (gc) at ($(g)!0.5!(c)$);
	\coordinate (cm) at ($(c)!0.6!(m)$);
	\coordinate (mi) at ($(m)!0.5!(i)$);
	\coordinate (ia) at ($(i)!0.5!(a)$);
	\coordinate (bf) at ($(b)!0.5!(f)$);
	\coordinate (fl) at ($(f)!0.6!(l)$);
	\coordinate (lp) at ($(l)!0.5!(p)$);
	\coordinate (ph) at ($(p)!0.5!(h)$);
	\coordinate (hd) at ($(h)!0.5!(d)$);
	\coordinate (do) at ($(d)!0.5!(o)$);
	\coordinate (ok) at ($(o)!0.5!(k)$);
	\coordinate (kb) at ($(k)!0.5!(b)$);
	\coordinate (lw) at ($(l)!0.5!(w)$);
	\coordinate (ws) at ($(w)!0.95!(s)$);
	\coordinate (sj) at ($(s)!0.5!(j)$);
	\coordinate (px) at ($(p)!0.5!(x)$);
	\coordinate (xt) at ($(x)!0.5!(t)$);
	\coordinate (tn) at ($(t)!0.5!(n)$);
	\coordinate (ku) at ($(k)!0.6!(u)$);
	\coordinate (uq) at ($(u)!0.5!(q)$);
	\coordinate (qi) at ($(q)!0.5!(i)$);
	\coordinate (ov) at ($(o)!0.5!(v)$);
	\coordinate (vr) at ($(v)!0.05!(r)$);
	\coordinate (rm) at ($(r)!0.5!(m)$);
	\coordinate (qa) at ($(q)!0.5!(a)$);
	\coordinate (bu) at ($(b)!0.5!(u)$);
	\coordinate (fw) at ($(f)!0.5!(w)$);
	\coordinate (es) at ($(e)!0.5!(s)$);
	\coordinate (dv) at ($(d)!0.5!(v)$);
	\coordinate (hx) at ($(h)!0.4!(x)$);
	\coordinate (gt) at ($(g)!0.5!(t)$);
	\coordinate (cr) at ($(c)!0.5!(r)$);

	\draw[decoration={markings, mark=at position 0.5 with {\arrow{>}}},postaction={decorate},color=black,ultra thick] (g)--(c);
	\draw[decoration={markings, mark=at position 0.5 with {\arrow{>}}},postaction={decorate},color=black,ultra thick,dashed] (d)--(o); 
	\draw[decoration={markings, mark=at position 0.5 with {\arrow{>}}},postaction={decorate},color=black,ultra thick] (s)--(e);
	
	\draw[decoration={markings, mark=at position 0.5 with {\arrow{>}}},postaction={decorate},color=red,ultra thick] (c)--(m); 
	\draw[decoration={markings, mark=at position 0.5 with {\arrow{>}}},postaction={decorate},color=red,ultra thick,dashed] (o)--(k);
	\draw[decoration={markings, mark=at position 0.5 with {\arrow{>}}},postaction={decorate},color=red,ultra thick] (l)--(w);
	
	\draw[decoration={markings, mark=at position 0.5 with {\arrow{>}}},postaction={decorate},color=green,ultra thick] (m)--(i);
	\draw[decoration={markings, mark=at position 0.5 with {\arrow{>}}},postaction={decorate},color=green,ultra thick,dashed] (k)--(b);
	\draw[decoration={markings, mark=at position 0.5 with {\arrow{>}}},postaction={decorate},color=green,ultra thick] (t)--(n);
	
	\draw[decoration={markings, mark=at position 0.5 with {\arrow{>}}},postaction={decorate},color=blue,ultra thick] (i)--(a);
	\draw[decoration={markings, mark=at position 0.5 with {\arrow{>}}},postaction={decorate},color=blue,ultra thick] (b)--(f);
	\draw[decoration={markings, mark=at position 0.5 with {\arrow{>}}},postaction={decorate},color=blue,ultra thick,dashed] (h)--(x);
	
	\draw[decoration={markings, mark=at position 0.5 with {\arrow{>}}},postaction={decorate},color=orange,ultra thick] (a)--(e);
	\draw[decoration={markings, mark=at position 0.5 with {\arrow{>}}},postaction={decorate},color=orange,ultra thick,dashed] (f)--(l);
	\draw[decoration={markings, mark=at position 0.5 with {\arrow{>}}},postaction={decorate},color=orange,ultra thick] (r)--(c);
	
	\draw[decoration={markings, mark=at position 0.5 with {\arrow{>}}},postaction={decorate},color=teal,ultra thick] (e)--(j);
	\draw[decoration={markings, mark=at position 0.5 with {\arrow{>}}},postaction={decorate},color=teal,ultra thick] (l)--(p);
	\draw[decoration={markings, mark=at position 0.5 with {\arrow{>}}},postaction={decorate},color=teal,ultra thick,dashed] (o)--(v);
	
	\draw[decoration={markings, mark=at position 0.5 with {\arrow{>}}},postaction={decorate},color=brown,ultra thick] (j)--(n);
	\draw[decoration={markings, mark=at position 0.5 with {\arrow{>}}},postaction={decorate},color=brown,ultra thick,dashed] (p)--(h);
	\draw[decoration={markings, mark=at position 0.5 with {\arrow{>}}},postaction={decorate},color=brown,ultra thick] (q)--(i);
	
	\draw[decoration={markings, mark=at position 0.5 with {\arrow{>}}},postaction={decorate},color=pink,ultra thick] (n)--(g);
	\draw[decoration={markings, mark=at position 0.5 with {\arrow{>}}},postaction={decorate},color=pink,ultra thick,dashed] (h)--(d);
	\draw[decoration={markings, mark=at position 0.5 with {\arrow{>}}},postaction={decorate},color=pink,ultra thick] (b)--(u);
	
	\draw[decoration={markings, mark=at position 0.5 with {\arrow{>}}},postaction={decorate},color=gray,ultra thick] (t)--(x);
	\draw[decoration={markings, mark=at position 0.5 with {\arrow{>}}},postaction={decorate},color=gray,ultra thick] (m)--(r);
	\draw[decoration={markings, mark=at position 0.5 with {\arrow{>}}},postaction={decorate},color=gray,ultra thick] (w)--(f);
	
	\draw[decoration={markings, mark=at position 0.5 with {\arrow{>}}},postaction={decorate},color=lightgray,ultra thick] (q)--(u);
	\draw[decoration={markings, mark=at position 0.5 with {\arrow{>}}},postaction={decorate},color=lightgray,ultra thick] (j)--(s);
	\draw[decoration={markings, mark=at position 0.5 with {\arrow{>}}},postaction={decorate},color=lightgray,ultra thick,dashed] (v)--(d);
	
	\draw[decoration={markings, mark=at position 0.5 with {\arrow{>}}},postaction={decorate},color=cyan,ultra thick] (g)--(t);
	\draw[decoration={markings, mark=at position 0.5 with {\arrow{>}}},postaction={decorate},color=cyan,ultra thick,dashed] (u)--(k);
	\draw[decoration={markings, mark=at position 0.5 with {\arrow{>}}},postaction={decorate},color=cyan,ultra thick] (s)--(w);
	
	\draw[decoration={markings, mark=at position 0.5 with {\arrow{>}}},postaction={decorate},color=lime,ultra thick] (a)--(q);
	\draw[decoration={markings, mark=at position 0.5 with {\arrow{>}}},postaction={decorate},color=lime,ultra thick] (x)--(p);
	\draw[decoration={markings, mark=at position 0.5 with {\arrow{>}}},postaction={decorate},color=lime,ultra thick,dashed] (r)--(v);

	\draw (ae) node[above right]{$1$};
	\draw (ej) node[below left]{$1$};
	\draw (jn) node[below right]{$1$};
	\draw (ng) node[below right]{$1$};
	\draw (gc) node[left]{$1$};
	\draw (cm) node[right]{$1$};
	\draw (mi) node[below]{$1$};
	\draw (ia) node[above left]{$1$};
	\draw (bf) node[right]{$a$};
	\draw (fl) node[left]{$a$};
	\draw (lp) node[above]{$a$};
	\draw (ph) node[below right]{$a$};
	\draw (hd) node[left]{$a$};
	\draw (do) node[above right]{$a$};
	\draw (ok) node[above]{$a$};
	\draw (kb) node[above left]{$a$};
	\draw (lw) node[above right]{$ab$};
	\draw (ws) node[right]{$c$};
	\draw (sj) node[above]{$b$};
	\draw (px) node[above]{$ba$};
	\draw (xt) node[above left]{$1$};
	\draw (tn) node[below left]{$b$};
	\draw (ku) node[above]{$ba$};
	\draw (uq) node[below right]{$1$};
	\draw (qi) node[below right]{$b$};
	\draw (ov) node[below]{$ab$};
	\draw (vr) node[left]{$c$};
	\draw (rm) node[below left]{$b$};
	\draw (qa) node[left]{$1$};
	\draw (bu) node[right]{$ba$};
	\draw (fw) node[above right]{$c$};
	\draw (es) node[right]{$c$};
	\draw (dv) node[left]{$c$};
	\draw (hx) node[above right]{$ba$};
	\draw (gt) node[left]{$1$};
	\draw (cr) node[below left]{$c$};
\end{tikzpicture}
\end{subfigure}
\caption{The 0-cells and 1-cells of $\mathcal{DV}_3$. Left: Representative 0-cells: black, red, cyan, green, brown, pink. Right: Representative 1-cells: black, red, green, blue, orange, teal, brown, pink, gray, lightgray, cyan, lime.}\label{fig7}
\end{figure}

\begin{figure}[h!]
\begin{subfigure}[h]{0.3\textwidth}
\centering
\begin{tikzpicture}[x={(0cm,-1cm)},y={(1cm,0cm)},z={(-3.85mm,-3.85mm)},scale=1.75]
	\coordinate (a) at (0.4142,1,1);
	\coordinate (b) at (0.4142,1,-1);
	\coordinate (c) at (0.4142,-1,1);
	\coordinate (d) at (0.4142,-1,-1);
	\coordinate (e) at (-0.4142,1,1);
	\coordinate (f) at (-0.4142,1,-1);
	\coordinate (g) at (-0.4142,-1,1);
	\coordinate (h) at (-0.4142,-1,-1);
	\coordinate (i) at (1,0.4142,1);
	\coordinate (j) at (-1,0.4142,1);
	\coordinate (k) at (1,0.4142,-1);
	\coordinate (l) at (-1,0.4142,-1);
	\coordinate (m) at (1,-0.4142,1);
	\coordinate (n) at (-1,-0.4142,1);
	\coordinate (o) at (1,-0.4142,-1);
	\coordinate (p) at (-1,-0.4142,-1);
	\coordinate (q) at (1,1,0.4142);
	\coordinate (r) at (1,-1,0.4142);
	\coordinate (s) at (-1,1,0.4142);
	\coordinate (t) at (-1,-1,0.4142);
	\coordinate (u) at (1,1,-0.4142);
	\coordinate (v) at (1,-1,-0.4142);
	\coordinate (w) at (-1,1,-0.4142);
	\coordinate (x) at (-1,-1,-0.4142);
	
	\coordinate (f21) at ($0.125*(a)+0.125*(i)+0.125*(m)+0.125*(c)+0.125*(g)+0.125*(n)+0.125*(j)+0.125*(e)$);
	\coordinate (f21a) at ($0.125*(b)+0.125*(f)+0.125*(l)+0.125*(p)+0.125*(h)+0.125*(d)+0.125*(o)+0.125*(k)$);
	\coordinate (f22) at ($0.125*(t)+0.125*(n)+0.125*(j)+0.125*(s)+0.125*(w)+0.125*(l)+0.125*(p)+0.125*(x)$);
	\coordinate (f22b) at ($0.125*(r)+0.125*(m)+0.125*(i)+0.125*(q)+0.125*(u)+0.125*(k)+0.125*(o)+0.125*(v)$);
	\coordinate (f23) at ($0.125*(a)+0.125*(q)+0.125*(u)+0.125*(b)+0.125*(f)+0.125*(w)+0.125*(s)+0.125*(e)$);
	\coordinate (f230) at ($0.125*(c)+0.125*(r)+0.125*(v)+0.125*(d)+0.125*(h)+0.125*(x)+0.125*(t)+0.125*(g)$);
	\coordinate (f24) at ($0.333*(g)+0.333*(n)+0.333*(t)$);
	\coordinate (f24ba) at ($0.333*(b)+0.333*(u)+0.333*(k)$);
	\coordinate (f25) at ($0.333*(c)+0.333*(m)+0.333*(r)$);
	\coordinate (f25ab) at ($0.333*(w)+0.333*(f)+0.333*(l)$);
	\coordinate (f26) at ($0.333*(a)+0.333*(i)+0.333*(q)$);
	\coordinate (f26ba) at ($0.333*(p)+0.333*(h)+0.333*(x)$);
	\coordinate (f27) at ($0.333*(s)+0.333*(e)+0.333*(j)$);
	\coordinate (f27ab) at ($0.333*(d)+0.333*(o)+0.333*(v)$);

	\draw (a)--(i)--(m)--(c)--(g)--(n)--(j)--(e)--(a);
	\draw (b)--(f);
	\draw[dashed,opacity=0.5] (f)--(l);
	\draw[dashed,opacity=0.5] (p)--(h)--(d)--(o);
	\draw[dashed,opacity=0.5] (k)--(b);
	\draw (a)--(e)--(s)--(w)--(f)--(b)--(u)--(q)--(a);
	\draw (c)--(g)--(t)--(x);
	\draw[dashed,opacity=0.5] (x)--(h)--(d)--(v)--(r);
	\draw (r)--(c);
	\draw (n)--(t)--(x)--(p)--(l)--(w)--(s)--(j);
	\draw (i)--(m)--(r);
	\draw[dashed,opacity=0.5] (r)--(v)--(o)--(k)--(u);
	\draw (u)--(q)--(i);

	\fill[fill=red,opacity=0.7] (a)--(i)--(m)--(c)--(g)--(n)--(j)--(e);
	\fill[fill=blue,opacity=0.7] (t)--(x)--(p)--(l)--(w)--(s)--(j)--(n);
	\fill[fill=green,opacity=0.7] (a)--(e)--(s)--(w)--(f)--(b)--(u)--(q);
	\fill[fill=brown,opacity=0.7] (n)--(t)--(g);
	\fill[fill=teal,opacity=0.7] (m)--(c)--(r);
	\fill[fill=orange,opacity=0.7] (a)--(i)--(q);
	\fill[fill=pink,opacity=0.7] (s)--(j)--(e);
\end{tikzpicture}
\end{subfigure}
\hspace{3cm}
\begin{subfigure}[h]{0.3\textwidth}
\centering
\begin{tikzpicture}[x={(0cm,-1cm)},y={(1cm,0cm)},z={(-3.85mm,-3.85mm)},scale=1.75]
	\coordinate (a) at (0.4142,1,1);
	\coordinate (b) at (0.4142,1,-1);
	\coordinate (c) at (0.4142,-1,1);
	\coordinate (d) at (0.4142,-1,-1);
	\coordinate (e) at (-0.4142,1,1);
	\coordinate (f) at (-0.4142,1,-1);
	\coordinate (g) at (-0.4142,-1,1);
	\coordinate (h) at (-0.4142,-1,-1);
	\coordinate (i) at (1,0.4142,1);
	\coordinate (j) at (-1,0.4142,1);
	\coordinate (k) at (1,0.4142,-1);
	\coordinate (l) at (-1,0.4142,-1);
	\coordinate (m) at (1,-0.4142,1);
	\coordinate (n) at (-1,-0.4142,1);
	\coordinate (o) at (1,-0.4142,-1);
	\coordinate (p) at (-1,-0.4142,-1);
	\coordinate (q) at (1,1,0.4142);
	\coordinate (r) at (1,-1,0.4142);
	\coordinate (s) at (-1,1,0.4142);
	\coordinate (t) at (-1,-1,0.4142);
	\coordinate (u) at (1,1,-0.4142);
	\coordinate (v) at (1,-1,-0.4142);
	\coordinate (w) at (-1,1,-0.4142);
	\coordinate (x) at (-1,-1,-0.4142);
	
	\coordinate (f21) at ($0.125*(a)+0.125*(i)+0.125*(m)+0.125*(c)+0.125*(g)+0.125*(n)+0.125*(j)+0.125*(e)$);
	\coordinate (f21a) at ($0.125*(b)+0.125*(f)+0.125*(l)+0.125*(p)+0.125*(h)+0.125*(d)+0.125*(o)+0.125*(k)$);
	\coordinate (f22) at ($0.125*(t)+0.125*(n)+0.125*(j)+0.125*(s)+0.125*(w)+0.125*(l)+0.125*(p)+0.125*(x)$);
	\coordinate (f22b) at ($0.125*(r)+0.125*(m)+0.125*(i)+0.125*(q)+0.125*(u)+0.125*(k)+0.125*(o)+0.125*(v)$);
	\coordinate (f23) at ($0.125*(a)+0.125*(q)+0.125*(u)+0.125*(b)+0.125*(f)+0.125*(w)+0.125*(s)+0.125*(e)$);
	\coordinate (f230) at ($0.125*(c)+0.125*(r)+0.125*(v)+0.125*(d)+0.125*(h)+0.125*(x)+0.125*(t)+0.125*(g)$);
	\coordinate (f24) at ($0.333*(g)+0.333*(n)+0.333*(t)$);
	\coordinate (f24ba) at ($0.333*(b)+0.333*(u)+0.333*(k)$);
	\coordinate (f25) at ($0.333*(c)+0.333*(m)+0.333*(r)$);
	\coordinate (f25ab) at ($0.333*(w)+0.333*(f)+0.333*(l)$);
	\coordinate (f26) at ($0.333*(a)+0.333*(i)+0.333*(q)$);
	\coordinate (f26ba) at ($0.333*(p)+0.333*(h)+0.333*(x)$);
	\coordinate (f27) at ($0.333*(s)+0.333*(e)+0.333*(j)$);
	\coordinate (f27ab) at ($0.333*(d)+0.333*(o)+0.333*(v)$);

	\draw (a)--(i)--(m)--(c)--(g)--(n)--(j)--(e)--(a);
	\draw (b)--(f);
	\draw[dashed,opacity=0.5] (f)--(l);
	\draw[dashed,opacity=0.5] (p)--(h)--(d)--(o);
	\draw[dashed,opacity=0.5] (k)--(b);
	\draw (a)--(e)--(s)--(w)--(f)--(b)--(u)--(q)--(a);
	\draw (c)--(g)--(t)--(x);
	\draw[dashed,opacity=0.5] (x)--(h)--(d)--(v)--(r);
	\draw (r)--(c);
	\draw (n)--(t)--(x)--(p)--(l)--(w)--(s)--(j);
	\draw (i)--(m)--(r);
	\draw[dashed,opacity=0.5] (r)--(v)--(o)--(k)--(u);
	\draw (u)--(q)--(i);
	
	\draw[pattern=north west lines,pattern color=red,opacity=0.5] (b)--(f)--(l)--(p)--(h)--(d)--(o)--(k);
	\draw[pattern=north west lines,pattern color=blue,opacity=0.5] (i)--(q)--(u)--(k)--(o)--(v)--(r)--(m);
	\draw[pattern=north west lines,pattern color=green,opacity=0.5] (c)--(g)--(t)--(x)--(h)--(d)--(v)--(r);
	\draw[pattern=north east lines,pattern color=brown,opacity=0.5] (b)--(u)--(k);
	\draw[pattern=north east lines,pattern color=teal,opacity=0.5] (w)--(l)--(f);
	\draw[pattern=north east lines,pattern color=orange,opacity=0.5] (p)--(x)--(h);
	\draw[pattern=north east lines,pattern color=pink,opacity=0.5] (d)--(o)--(v);
	
	\draw (f21a) node{$a$};
	\draw (f22b) node{$b$};
	\draw (f230) node{$c$};
	\draw (f24ba) node{$ba$};
	\draw (f25ab) node[above right]{$ab$};
	\draw (f26ba) node{$ba$};
	\draw (f27ab) node{$ab$};
\end{tikzpicture}
\end{subfigure}
\caption{The 2-cells of $\mathcal{DV}_3$. Left: Representative 2-cells: red, blue, green, brown, teal, orange, pink. Right: Some translates of the representing 2-cells.}\label{fig8}
\end{figure}
\end{appendix}

\subsection*{Acknowledgments}
We sincerely thank the referee for carefully reading of the manuscript and for the useful remarks, that helped to greatly improve the paper. Moreover, the referee's suggestion of looking at the case of lens spaces was more than welcome, and led us to the enlightening Example \ref{lens_spaces}, so we are much grateful for this as well.

\printbibliography

\end{document}